\documentclass[reqno]{amsart}

\usepackage{general,strings,txfonts}
\usepackage[bookmarks, colorlinks=true, pdfstartview=FitV, linkcolor=blue, citecolor=blue, urlcolor=blue]{hyperref}

\newcommand{\tubezf}{tubular\xspace}
\newcommand{\arxiv}[1]{\href{http://arxiv.org/abs/#1}{\texttt{arXiv:#1}}}
\newcommand{\abscissa}{\ensuremath{D}\xspace}

\numberwithin{equation}{section} \numberwithin{theorem}{section}

\title[Tube formulas and complex dimensions]{Tube formulas and complex dimensions \linebreak of self-similar tilings}

\author{Michel L. Lapidus}
\address{Department of Mathematics, University of California, Riverside, CA 92521-0135 USA}
\email{lapidus@math.ucr.edu}
\urladdr{http://www.math.ucr.edu/$\sim$lapidus}

\author{Erin P. J. Pearse}
\address{Department of Mathematics, Cornell University, Ithaca, NY 14853-4201 USA}
\email{erin@math.cornell.edu}
\urladdr{http://www.math.uiowa.edu/$\sim$epearse/}

\curraddr{University of Iowa\\Iowa City, IA 52242-1419 USA}
\email{erin-pearse@uiowa.edu}
   
\thanks{The work of MLL was partially supported by the US National
Science Foundation under the research grants DMS-0070497 and
DMS-0707524. The work of EPJP was partially supported by the University of Iowa Department of Mathematics NSF VIGRE grant DMS-0602242.}

\keywords{Complex dimensions, zeta functions, tube formula, fractal Steiner formula, inradius, self-similar tiling, curvature matrix, generating function for the geometry, distributional explicit formula, fractal string.}

\subjclass{11M41, 28A12, 28A75, 28A80, 52A39, 52C07 (primary), 11M36, 28A78, 28D20, 42A16, 42A75, 52A20, 52A38 (secondary)}

\allowdisplaybreaks

\begin{document}

\begin{abstract}
We use the self-similar tilings constructed in \cite{SST} to define a generating function for the geometry of a self-similar set in Euclidean space. This \tubezf zeta function encodes scaling and curvature properties related to the complement of the fractal set, and the associated system of mappings. This allows one to obtain the complex dimensions of the self-similar tiling as the poles of the \tubezf zeta function and hence develop a tube formula for self-similar tilings in \bRd. The resulting power series in \ge is a fractal extension of Steiner's classical tube formula for convex bodies $K \ci \bRd$. Our sum has coefficients related to the curvatures of the tiling, and contains terms for each integer $i=0,1,\dots,d-1$, just as Steiner's does. However, our formula also contains a term for each complex dimension. This provides further justification for the term ``complex dimension''. It also extends several aspects of the theory of fractal strings to higher dimensions and sheds new light on the tube formula for fractals strings obtained in \cite{FGCD}.
\end{abstract}

\maketitle
\setcounter{tocdepth}{1}\tableofcontents


\section{Introduction}
\label{ssec:overview}

The main result of this paper is a tube formula for a certain class of fractal sets. Here, a tube formula for a bounded set $A \ci \bRd$ is an explicit expression for the $d$-dimensional volume of the inner \ge-neighbourhood of $A$, i.e.,
\begin{align}
  \label{def:tube-formula}
  V(A,\ge) = \vol[d]\{x \in A \suth dist(x,\del A) \leq \ge\}.
\end{align}
Such formulas have myriad applications to geometry and have roots in the results of Steiner (when $A$ is convex) and Weyl (when $A$ is a smooth manifold). 
In order to explain how our result is an extension of these classical formulas, and how it is related to the development of a notion of curvature for fractal sets, we give a brief encapsulation of Steiner's theorem. Here, the Minkowski sum of two subsets $A$ and $B$ of \bRd is denoted by $A+B = \{x \in \bRd \suth x=a+b \text{ for } a \in A, b \in B\}.$ Then, if $B^k$ is the closed $k$-dimensional unit ball in $\bR^k$, one can denote the \emph{outer} \ge-neighbourhood of a set $A \ci \bRd$ by
\begin{align*}
  A + \ge B^d = \{x \in \bRd \suth dist(x,A) \leq \ge\}.
\end{align*}

\begin{theorem}[Steiner's formula]
  \label{thm:Steiners-formula-basic}
  If 
  $A \ci \bRd$ is convex and compact, then the $d$-dimensional volume of $A+\ge B^d$ is given by
  \begin{align}\label{eqn:Steiner-formula-with-inv-meas}
    \vol[d](A+\ge B^d) = \sum_{k=0}^{d} \gm_k(A) \vol[d-k](B^{d-k}) \ge^{d-k},
  \end{align}
  where $\gm_k$ is the renormalized $k$-dimensional \emph{intrinsic volume}.
\end{theorem}

Note that this formula is simply a polynomial in \ge; the coefficients are constants determined by the curvature of $A$ (and the unit ball). Up to some normalizing constant, the $k$-dimensional intrinsic volume is the same thing as the \nth[(d-k)] Quermassintegral (from Minkowski's theory of mixed volumes). The valuation $\gm_k$ can be defined via integral geometry as the average measure of orthogonal projections to $(d-k)$-dimensional subspaces; see \cite[Chap.~7]{Rota}. For now, we note that (up to a multiplicative constant), there is a correspondence
\begin{center}
\begin{tabular}{rrlrrl}
  $\gm_0$ & $\sim$ &Euler characteristic, \qq \qq & $\gm_{d-1}$ & $\sim$ &surface area, \\
  $\gm_1$ & $\sim$ &mean width, & $\gm_d$ & $\sim$ &volume,
\end{tabular}
\end{center}
see \cite[\S4.2]{Schn2} for further details.

When $A$ is sufficiently regular (i.e., when its boundary is a $C^2$ surface), these coefficients can be given in terms of curvature tensors, and in fact Steiner's tube formula coincides with the one obtained by Weyl in \cite{We}. In \cite{Fed}, Federer unified the tube formulas of Steiner (for convex bodies, as described in \cite[Ch.~4]{Schn2}) and of Weyl (for smooth submanifolds, as described in \cite{Gr} and \cite{We}) and extended these results to sets of \emph{positive reach}.\footnote{A set $A$ has positive reach iff there is some $\gd>0$ such that any point $x$ within \gd of $A$ has a unique metric projection to $A$, i.e., that there is a unique point $A$ minimizing $dist(x,A)$. Equivalently, every point $q$ on the boundary of $A$ lies on a sphere of radius \gd which intersects $\del A$ only at $q$.}
It is worth noting that Weyl's tube formula for smooth submanifolds of \bRd is expressed as a polynomial in \ge with coefficients defined in terms of curvatures (in the classical sense) that are \emph{intrinsic} to the submanifold \cite{We}. See \cite[\S6.6--6.9]{BeGo} and the book \cite{Gr}.
Federer's tube formula has since been extended in various directions by a number of researchers in integral geometry and geometric measure theory, including \cite{Schn1}, \cite{Schn2}, \cite{Za1}, \cite{Za2}, \cite{Fu1},\cite{Fu2}, \cite{Stacho}, and most recently (and most generally) in \cite{HuLaWe}. The books \cite{Gr} and \cite{Schn2} contain extensive endnotes with further information and many other references.

\pgap

Note that \eqref{eqn:Steiner-formula-with-inv-meas} gives the volume of the set of points which are within \ge of $A$, including the points of $A$. If we denote the \emph{exterior \ge-neighbourhood} of $A$ by
\begin{align} \label{eqn:ext-nbd-of-A}
  A_\ge^{ext} := \{x \in \bRd \less A \suth dist(x,A) \leq \ge\},
\end{align}
then it is immediately clear that omitting the \nth[d] term gives
\begin{align}
  \vol[d](A_\ge^{ext}) = \sum_{k=0}^{d-1} C_k(A) \ge^{d-k}
      \label{def:Steiner ext formula with inv meas}
\end{align}
with $C_k(A) = \gm_k(A) \vol[d-k](B^{d-k})$. These coefficients $C_k(A)$ are called the \emph{total curvatures} of $A$ and they are key geometric invariants. Even more importantly, $C_k(A)$ can be localized and understood as the \emph{curvature measures} described in \cite{Fed} and \cite[Ch.~4]{Schn2}. In this case, for a Borel set $\gb \ci \bRd$, one has
\begin{align}
  \vol[d]\{x \in A_\ge^{ext} \suth p(x,A) \in \gb\}
  = \sum_{k=0}^{d-1} C_k(A,\gb) \ge^{d-k}
    \label{def:Steiner ext formula localized}
\end{align}
where $p(x,A)$ is the metric projection of $x$ to $A$, that is, the closest point of $A$ to $x$.
In fact, the curvature measures are obtained axiomatically in \cite{Schn2} as the coefficients of the tube formula, and it is this approach that we hope to emulate in a forthcoming work, based on a localized version of the tube formula obtained in the present paper. In other words, we believe that $\crv_k$ (introduced in \eqref{def:preview-c_wj} below) may also be understood as a (total) curvature, in a suitable sense, and we expect that $\crv_k$ can be localized as a curvature measure (or rather, current). A more rigorous formulation of these ideas is currently underway in \cite{FCM}.

\pgap

\subsubsection*{Our main result}
Since it is somewhat arduous to express our tube formula with complete precision, we present a somewhat simplified version of it here, so that the key features are not overshadowed by technical details (the exact hypotheses necessary for this formulation are given in Corollary~\ref{cor:self-similar-one-generator}).

\begin{theorem}
  \label{thm:preview-of-simplified-tube-formula}
  The $d$-dimensional volume of the tubular neighbourhood of a sufficiently nice self-similar fractal set \attr is given by the following distributional explicit formula:
  \begin{align}
    \label{eqn:preview-of-simplified-tube-formula}
    V(\attr,\ge) &= \negsp[10]\sum_{\gw \in \Ds \cup \intDim} \negsp[10] c_\gw \ge^{d-\gw},
  \end{align}
  where for each fixed $\gw \in \Ds$,
  \begin{align}
    \label{def:preview-c_wj}
    c_\gw
    &:=  \res{\gzs(s)} \sum_{k=0}^d \frac{\genir^{\gw-k}}{\gw-k} \crv_k.
  \end{align}
\end{theorem}
  Here, \gzs is a zeta function which encodes scaling information about \attr, \Ds is the set of poles of \gzs (the \emph{complex dimensions}\footnote{Actually, elsewhere in this paper, \emph{complex dimensions} refers to the set $\DT = \Ds \cup \intDim$ (the poles of \gzT) and we refer to \Ds as the \emph{scaling} complex dimensions; see Definition~\ref{def:visible complex dimensions of a fractal spray}.}), 
  and \genir and $\crv_k$ are geometric data obtained from the complement of \attr. Roughly speaking, \genir is the size of the largest bounded component of $\bRd \less \attr$ (the \emph{generator}), and $\crv_k$ are ``curvatures'' of the same set, in some sense.  Also, ``sufficiently nice'' requires that the open set condition is satisfied for a set whose boundary is contained in \attr; cf.~\cite[Thm.~6.2]{GeometryOfSST}.
  
  Theorem~\ref{thm:preview-of-simplified-tube-formula} is a special case of a more fundamental result which is valid under much more general hypotheses: the tube formula for fractal sprays given in Theorem~\ref{thm:fractal-spray-tube-formula}:

\begin{theorem}
  \label{thm:preview-of-tube-formula}
  The $d$-dimensional volume of the inner tubular neighbourhood of a fractal spray \tiling with a generator \gen is given by the following distributional explicit formula:
  \begin{align}
    \label{eqn:preview-of-tube-formula}
    V(\tiling,\ge) &= \sum_{\gw \in \DT} \negsp[1] \res{\gzT(\ge,s)},
  \end{align}
  where 
  \begin{align}\label{eqn:preview-tube-zeta}
    \gzT(\ge,s) := \ge^{d-s} \gzs(s) \sum_{k=0}^d \frac{\genir^{s-k}}{s-k} \crv_{k} 
  \end{align}
  and $\DT := \Ds \cup \intDim$ is the set of poles of the \tubezf zeta function \gzT. 
\end{theorem}
In \eqref{eqn:preview-tube-zeta}, the numbers $\crv_0, \crv_1, \dots, \crv_d$ refer to a chosen ``Steiner-like representation'' of \gen (an expression of the inner tube for \gen which satisfies certain very mild conditions discussed in Definition~\ref{def:Steiner-like}.
 
Our tube formula extends previous results in two ways. On one hand, it provides a fractal analogue of the classical Steiner formula of convex geometry. 
On the other hand, the tube formula \eqref{eqn:preview-of-tube-formula} also provides a natural higher-dimensional analogue of the tube formula for fractal strings obtained in \cite{FGCD} and recalled in \eqref{eqn:1-dim tube formula}. The present work can be considered as a further step towards a higher-dimensional theory of fractal strings and their complex dimensions, especially in the self-similar case, following upon \cite[\S10.2 and \S10.3]{FGNT1}, and our earlier paper \cite{KTF}. This is discussed further in \S\ref{ssec:The tube formula for fractal strings} and in Remark~\ref{rem:comparison of Koch tiling tube to Koch tube}.

To emphasize the present analogy with \eqref{eqn:preview-of-simplified-tube-formula}, consider that, with the obvious change of notation, Steiner's formula \eqref{def:Steiner ext formula with inv meas} may be rewritten
\begin{align}
  \label{eqn:steiner_conceptual}
  \vol[d](A_\ge^{ext}) = \sum_{k \in \{0,1,\dots,d-1\}} c_k \ge^{d-k}.
\end{align}
The obvious similarities between the tube formulas \eqref{eqn:preview-of-simplified-tube-formula} and \eqref{eqn:steiner_conceptual} is striking. Our tube formula is a ``fractal power series'' in \ge, rather than just a polynomial in \ge as in Steiner's (and Weyl's) formula. Moreover, our series is summed not just over the `integral dimensions' \intDim, but also over the countable set \Ds of complex dimensions. The coefficients $c_\gw$ of the tube formula are expressed in terms of the `curvatures' and the inradii of the generators of the tiling. It is intriguing to consider that the extra terms appearing in our formula are oscillatory, as evinced by the purely imaginary components of the complex dimensions. Moreover, fractals are objects with (multiplicative) geometric oscillations, in the sense that the same geometric patterns repeat, at different scales of magnification. This is the essential theme of the theory of complex dimensions as expounded in \cite{FGCD}.

The primary object of study in \cite{FGCD} is a \emph{fractal string}, a countable collection $L = \{L_j\}_{j=1}^\iy$ of disjoint open intervals which form a bounded open subset of \bR. Due to the trivial geometry of such intervals, this reduces to studying the lengths of these intervals \mbox{$\sL = \{\ell_j\}_{j=1}^\iy$}, \label{page:defn-of-Ln} and the sequence \sL is also referred to as a fractal string. The tube formula for a fractal string \sL (and in particular, for a self-similar tiling in $\bR$) is defined to be $V(\sL,\ge) := V_L(\ge)$ and is shown to be essentially given by a sum of the form
\begin{align}
  \label{eqn:form of 1-dim tube formula}
  V(\sL,\ge) = \sum_{\gw \in \D_\sL \cup\{0\}} \negsp[10]  c_\gw \ge^{1-\gw}
\end{align}
in \cite{FGCD}, Thm.~8.1. Here, the sum is taken over the set of complex dimensions $\D_\sL = \{\text{poles of }\gzL\}$, and $c_\gw$ is given in terms of the residue of $\gzL(s)$ at $s=\gw$, the geometric zeta function of \sL (defined as the meromorphic continuation of the Dirichlet series $\sum_{j=1}^\iy \ell_j^s$, for $s \in \bC$). The definition $V(\sL,\ge) := V(L,\ge)$ is justified because, as is shown in \cite{LaPo1}, $V(L,\cdot)$ depends exclusively on \sL.

In \cite[\S1.4]{FGCD} (following \cite{LaPo2}), a \emph{fractal spray} is defined to be given by a nonempty bounded open set $\gen \ci \bRd$ called the \emph{generator} (or ``\emph{basic shape}'' in \cite{FGCD}), scaled by a fractal string \mbox{$\sL = \{\ell_j\}_{j=1}^\iy$}. That is, a fractal spray is a bounded open subset of \bRd which is the disjoint union of open sets $\gW_j$ for $j=1,2,\dots$, where $\gW_j$ is congruent to $\ell_j \gen$ (the homothetic of \gW by $\ell_j$) for each $\ell_j$. Thus, a fractal string is a fractal spray on the generator $\gen=(0,1)$, the unit interval. In the context of the current paper, a self-similar tiling is a union of fractal sprays on the generators $\gen_1,\dots,\gen_Q$, each scaled by a fixed \emph{self-similar} string. In fact, we first prove Theorem~\ref{thm:preview-of-tube-formula} for the more general case of fractal sprays, and then refine it to obtain the formula for self-similar tilings.

For fractal sets, the tube formula may also be used to assess the Minkowski measurability and to determine the Minkowski dimension (also called the \emph{box} or \emph{box-counting} dimension) and the value of the Minkowski content (when it exists); see Corollary~\ref{cor:meas and dichotomy} and Remark~\ref{rem:dimensions-and-measurability}. \emph{Minkowski dimension} of the boundary $\del A$ of a bounded open subset $A \ci \bRd$ is a real-valued extension of the usual (topological) notion of dimension defined by 
\begin{align*}
  D = \dim_\sM(\del A) := \inf \{\ga \geq 0 \suth V(A,\ge) = O(\ge^{1-\ga}) \text{ as } \ge \to 0^+\}.
\end{align*}
Minkowski dimension frequently coincides with Hausdorff dimension (for example, for fractal sets defined in terms of an iterated function system with mappings that do not overlap too much) and in general, $\dim_\sH(X) \leq \dim_\sM(X)$. The set $\del A$ is said to be \emph{Minkowski measurable}, with \emph{Minkowski content}
\begin{align*}
  \sM(\del A) := \lim_{\ge \to 0^+} V(A,\ge) \ge^{-(d-\dim_\sM(\del A))},
\end{align*}
if this limit exists and is both finite and strictly positive. 

Further (potential) applications and extensions of our results are discussed in \S\ref{sec:Concluding Remarks} and elsewhere in the paper; in particular, please see Remark~\ref{rem:combine-FGCD-with-this}. Some of the main results of this paper were announced in \cite{TFSST}.

\subsubsection*{Strategy for obtaining the tube formula} 
For simplicity, we limit ourselves here to the self-similar case.
In \cite{SST}, the second author has shown that a certain (self-similar) tiling \tiling is canonically associated with any self-similar set. This tiling is defined via the finite collection $\simt=\{\simt_j\}_{j=1}^J$ of contractive similarity transformations which defines the self-similar set. This tiling \tiling is essentially a decomposition of the complement of the unique self-similar set associated with \gF, and is reviewed in greater detail in \S\ref{sec:The self-similar tiling}. 
Several conditions are given in \cite{GeometryOfSST} which describe precisely when the tube formula of a self-similar set may be obtained from a tube formula for the self-similar tiling \tiling. The results of the present study allow one to obtain a tube formula for the (large) class of fractal sets satisfying the ``compatibility conditions'' of \cite[Thm.~4.4 or Thm.~6.2]{GeometryOfSST}.

At the heart of this paper is the \tubezf zeta function $\gzT(\ge,s)$ of a self-similar tiling \tiling. It will take some work before we are able to describe this meromorphic distribution-valued function precisely in \S\ref{sec:The tube formula for self-similar tilings}. The function \gzT is a generating function for the geometry of a self-similar tiling: it encodes the density of geometric states of a tiling, including curvature and scaling properties. The poles of \gzT are the \emph{complex dimensions} \DT of the tiling, and we obtain a tube formula for \tiling given as a sum over \DT of the residues of \gzT, taken at the complex dimensions. The complex dimensions generalize the Minkowski dimension in the sense that, for any tiling \tiling, $\sup\{\DT \cap \bR\}$ is equal to the Minkowski dimension of the underlying self-similar set \attr.

The first ingredient of \gzT is a \emph{scaling zeta function} $\gzs(s)$ which encodes the scaling properties of the tiling and is discussed in \S\ref{ssec:scaling-zet-fn}. This comparatively simple zeta function is the Mellin transform of a discrete \emph{scaling measure} \ghs which encodes the combinatorics of the scaling ratios of a self-similar tiling. More precisely, if one considers a composition of similarity mappings $\simt_j$, each with scaling ratio $r_j$, for $j=1,\dots,J$, then
\[\simt_w = \simt_{w_1w_2\dots w_n} = \simt_{w_n} \comp \dots \comp \simt_{w_2} \comp \simt_{w_1}\]
has scaling ratio $r_w = r_{w_1}r_{w_2} \dots r_{w_n}$, where $w_i \in \{1,2,\dots,J\}$. The measure \ghs is a sum of Dirac masses, where each mass is located at a reciprocal scaling ratio $r_w^{-1}$. The total mass of any point in the support of \ghs corresponds to the multiplicity with which such a scaling ratio can occur. The scaling zeta function \gzs is formally identical to the zeta functions studied in \cite{FGCD}. The function \gzs also allows us to define the scaling complex dimensions of a self-similar set in \bRd (as the poles of \gzs), and we find these dimensions to have the same structure as in the 1-dimensional case; see \S\ref{ssec:comparison to FGNT} and Remark~\ref{rem:dimensions-and-measurability}. The definition and properties of the scaling measure \ghs and zeta function \gzs is the subject of \S\ref{sec:The geometric zeta function of a self-similar tiling}.

The next ingredient of \gzT is a \emph{generator tube formula} \gtf. In \cite{SST}, it is shown that certain tiles $\gen_1,\dots,\gen_Q$ of \tiling are generators in the sense that any tile $\tile_n$ of \tiling is the image of some $\gen_q$ under some composition of the mappings $\simt_j$, i.e.,
\begin{align*}
  \tile_n \in \tiling \q\implies\q \tile_n = \simt_w(\gen_q),
\end{align*}
for some $\gen_q$ and some $w = w_1 w_2 \dots w_m$. In \S\ref{sec:tube formula for tilings}, we discuss the role of the generators and introduce the function \gtf which gives the inner tube formula for a generator in the sense of \eqref{def:tube-formula}. Moreover, appropriately parameterizing \gtf yields the inner tube formula for a scaled generator. Therefore, by integrating \gtf against \ghs, one obtains the total contribution of $\gen_q$ (and its images under the maps $\simt_w$) to the final tube formula $V_\tiling$. This is elaborated upon in \S\ref{ssec:Tilings with one generator}.

At last, the \tubezf zeta function of the tiling \gzT is assembled from the scaling zeta function, the tiling, and the terms appearing in \gtf. In some precise sense, \gzT is a generating function for the geometry of the self-similar tiling. Using \gzT, and following the distributional techniques and explicit formulas of \cite{FGCD}, we are able to obtain an explicit distributional tube formula for self-similar tilings.

\subsubsection*{Outline}
The rest of this paper is organized as follows.
\S\ref{sec:The self-similar tiling} contains a quick overview of the background material concerning self-similar tilings. \S\ref{ssec:The inradius} discusses how the notion of inradius describes the different scales of the tiling. \S\ref{ssec:Measures and zeta functions} defines the scaling measure, the scaling zeta function, and complex dimensions of a self-similar tiling. \S\ref{sec:tube formula for tilings} develops the tube formula for the generators of a self-similar tiling, and establishes the general form of $V(\tiling,\ge)$ in terms of this. \S\ref{sec:The distributional theory of fractal strings} reviews the explicit formulas for fractal strings which will be used in the proof of the main results. \S\ref{sec:The tube formula for self-similar tilings} defines the \tubezf zeta function of the tiling, and states and proves the tube formula for fractal sprays (a generalization of a tiling) given in Theorem~\ref{thm:fractal-spray-tube-formula}, from which the tube formula for self-similar tilings follows readily, and \S\ref{sec:Tube formula examples} discusses several examples illustrating the theory.
 Appendix \ref{app:Rigourous Definition of gzT} verifies the validity of the definition of the \tubezf zeta function \gzT. Finally, Appendix \ref{app:The Error Term and Estimate} verifies the distributional error term and its estimate, from Theorem~\ref{thm:fractal-spray-tube-formula}.

\begin{remark}[A note on the references]
  \label{rem:explaining the references}
  The primary references for this paper are \cite{SST} and the research monograph ``Fractal Geometry, Complex Dimensions and Zeta Functions: Geometry and spectra of fractal strings'' by Lapidus and van Frankenhuijsen \cite{FGCD}. This volume is essentially a revised and much expanded version of
  \cite{FGNT1}, by the same authors. The present paper cites \cite{FGCD} almost exclusively, so we provide the following partial correspondence between chapters for the aid of the reader:
  \begin{center}
    \begin{tabular}{l|l|l|l|l|l|}
      \cite{FGNT1} & Ch.~2 & Ch.~3& Ch.~4 & Ch.~6 & Ch.~10 \\ \hline
      \cite{FGCD} & Ch.~2--3 & Ch.~4 & Ch.~5 & Ch.~8 & Ch.~12
    \end{tabular}
  \end{center}
\end{remark}

\begin{remark}
  \label{rem:imaginary i}
  Throughout, we reserve the symbol $\ii = \sqrt{-1}$ for the
  imaginary number.
\end{remark}

\subsection{Acknowledgements}
\label{sec:acknowledgements}

The authors wish to thank Martina Z\"{a}hle for several helpful discussions on geometric measure theory and for bringing the reference \cite{HuLaWe} to our attention. Additionally, the authors would like to thank Steffen Winter for many helpful discussions and suggestions, and for finding mistakes in an earlier version of this paper. Steffen also suggested the term ``monophase'' (the authors had originally used the term ``diphase'' in an earlier draft of this paper).


\section{The Self-Similar Tiling}
\label{sec:The self-similar tiling}

This section provides an overview of the necessary background
material concerning self-similar tilings. Further details may be
found in \cite{SST}.

\begin{defn}
  A \emph{self-similar system} is a
  family $\{\simt_\j\}_{\j=1}^J$ (with $J \geq 2$) of
  contraction similitudes
  \[\simt_\j(x) := r_\j M_\j x + a_\j, \q \j=1,\dots,J.\]
  For $\j=1,\dots,J$, we have $0 < r_\j < 1, a_\j \in \bRd$, and $M_\j
  \in O(d)$, the orthogonal group of rigid rotations in $d$-dimensional
  Euclidean space \bRd. The number $r_\j$ is the \emph{scaling ratio}
  of $\simt_j$. For convenience, assume that
  \begin{equation}
    \label{eqn:scaling ratio ordering}
    1> r_1 \geq r_2 \geq \dots \geq r_J >0.
  \end{equation}
\end{defn}

It is well known that there is a unique nonempty compact subset
$\attr \ci \bRd$ satisfying the fixed-point equation
  \begin{equation}
    \label{eqn:fixed-pt_eqn}
    \attr = \simt(\attr) := \bigcup_{\j=1}^J \simt_\j(\attr).
  \end{equation}
This (self-similar) set $F$ is called the \emph{attractor} of \simt.
We abuse notation and let \simt denote both an operator on compacta
(as in \eqref{eqn:fixed-pt_eqn}) and the family $\{\simt_\j\}$.
Different self-similar systems may give rise to the same
self-similar set; therefore we emphasize the self-similar system and
its corresponding dynamics.

\begin{figure}
  \centering
  \scalebox{0.90}{\includegraphics{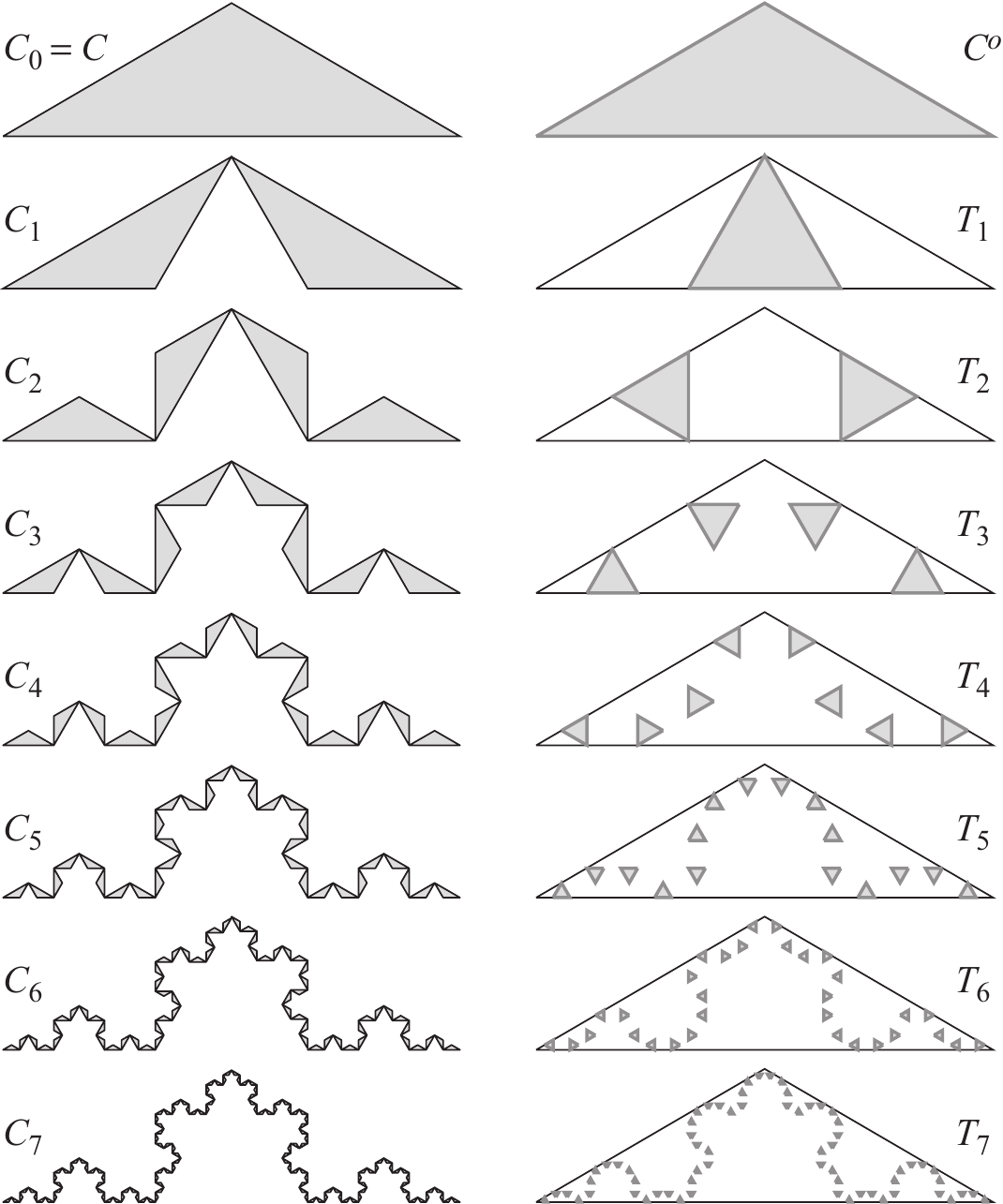}}
  \caption[Construction of the Koch tiling \sK.]
    {\captionsize Construction of the Koch tiling \sK.
    This example is discussed further in \S\ref{sec:Tube formula examples}.
    The tiling \sK has the single generator $\tileset_1 = \gen_1$, an
    equilateral triangle.}
  \label{fig:koch-contractions}
\end{figure}

It is shown in \cite{SST} that for a self-similar system satisfying the
tileset condition (see Definition~\ref{defn:tileset-condition}), there
exists a natural decomposition of $\hull \less \attr$ which is
produced by the system $\simt$. The construction of this tiling is
illustrated for a well-known example, the Koch curve, in
Figure~\ref{fig:koch-contractions}. It may help the reader to look at
this example before diving into the next paragraph and the thicket
of definitions therein. Further examples are depicted in
\S\ref{sec:Tube formula examples}.

Let $\hull:=[\attr]$ be the convex hull of \attr, and let $\tileset := \relint \hull$ be the relative interior of \hull. Iterates of the hull $\hull$ under $\simt$ are denoted
\begin{equation}
  \label{def:hullk}
  \hull_k := \simt^k(\hull) = \bigcup_{w \in \Wds_k} \simt_w(\hull),
\end{equation}
where $w=w_1\dots w_k$ is a \emph{word} in $\Wds_k := \{1,2,\dots,J\}^k$ and $\simt_w := \simt_{w_k} \comp \dots \comp \simt_{w_2} \comp \simt_{w_1}$. For future reference, let $\Wds := \bigcup_{k=1}^\iy \Wds_k$ be the set of all \emph{finite} words $w$ over the alphabet $\{1,2,\dots,J\}$.

\begin{defn}
  \label{defn:tileset-condition}
  The system satisfies the \emph{tileset condition} iff $\tileset
  \nsubseteq \simt(\hull)$ and
  \begin{equation}
    \label{eqn:tileset condition}
      \inter \simt_j(\hull) \cap \inter \simt_\ell(\hull) = \es,
      \qq j \neq \ell.
  \end{equation}
\end{defn}
This is a restriction on the overlap of the images of the mappings and implies (but is not equivalent to) the open set condition. For any system satisfying the tileset condition,
\begin{equation}
  \label{def:tileset1}
  \tileset_1 := \tileset \less \hull_1 
\end{equation}
is well-defined and nonempty, and hence so is $\tileset_k :=
\simt^k(\tileset_1)$. As an open set, $\tileset_1$ is a
disjoint union of connected open sets:
\begin{equation}
  \label{def:tileset1 as union}
  \tileset_1 = \gen_1 \cup \gen_2 \cup \dots \cup \gen_Q,
  \q \gen_p \cap \gen_q = \es, p \neq q.
\end{equation}

\begin{defn}
  \label{def:generators}
  The \emph{generators} of the tiling are the connected components
  of $\tileset_1$, i.e., the disjoint open sets $\{\gen_q\}$ in
  \eqref{def:tileset1 as union}.
\end{defn}
The number $Q$ of generators depends on the system \simt, not just
on \attr. In general, the number of connected components of an open
subset of \bRd may be countable; however, in this paper we assume
$Q<\iy$.

\begin{defn}
  \label{def:ss tiling}
  The \emph{self-similar tiling} of \simt is 
  \begin{equation}
    \label{eqn:def_tiles}
    \tiling = \{\tile_n\}_{n=1}^\iy = \{\simt_w(\gen_q) \suth w \in \Wds, q=1,\dots,Q\}.
  \end{equation}
\end{defn}
In \eqref{eqn:def_tiles}, the sequence $\{\tile_n\}$ is an
enumeration of the sets $\{\simt_w(\gen_q)\}$, and $\simt_w$ is
as in \eqref{def:hullk}.
We say \tiling is a \emph{tiling of $\hull \less \attr$} because the tiles
$\tile_n$ have disjoint interiors and \attr does not intersect the
interior of any $\tile_n$ (see Figure~\ref{fig:koch-tiled-in}):
\begin{align*}
  \hull = \bigcup\nolimits_n \cj{\tile_n}, \;\;
  \attr \ci \bigcup\nolimits_{n} \del \tile_n, \;\;
  \text{and} \;\;
  \tile_{n_1} \cap \tile_{n_2} = \del \tile_{n_1} \cap \del \tile_{n_2}.
\end{align*}


\section{The Inradius} \label{ssec:The inradius}

As alluded to in \eqref{def:tube-formula}, we are interested in that portion of a set which lies within \ge of its boundary.
\begin{defn}
  \label{def:V of a set}
  Given $\ge > 0$, the \emph{inner \ge-neighbourhood} of a bounded set $A \ci \bRd$,
  $d \geq 1$, is
  \begin{equation}\label{eqn:def:A_ge}
    A_\ge := \{x \in A \suth dist(x,\del A) \leq \ge\},
  \end{equation}
  where $\del A$ is the boundary of $A$. We are primarily interested in
  the $d$-dimensional Lebesgue measure of $A_\ge$, denoted
  $V_{A}(\ge) := \vol[d](A_\ge)$.
\end{defn}

\begin{remark}\label{rem:inner-instead-of-outer}
  The primary reason we have worked with the inner \ge-neighbourhood instead of the exterior (as in \eqref{eqn:ext-nbd-of-A}) is that it is more intrinsic to the set; it makes the computation independent of the embedding of \tiling into \bRd. At least, this should be the case, provided the `curvature' terms $\crv_k$ of Definition~\ref{def:monophase} are also intrinsic. As a practical bonus, working with the inner \ge-neighbourhood allows us to avoid potential issues with the intersections of the \ge-neighbourhoods of different components.
\end{remark}

It is clear that if $A$ is a bounded set, $A \ci A_\ge$ for sufficiently large \ge. Alternatively, it is apparent that for a fixed $\ge>0$, any sufficiently small set will be entirely contained within its \ge-neighbourhood. The notion of \emph{inradius} allows us to see when this phenomenon occurs.


\begin{defn}
  The \emph{inradius} \gr of a set $A$ is
  \begin{equation}
    \label{def:inradius}
    \gr = \gr(A)
      := \sup \{\ge > 0 \suth \exists x \text{ with } B(x,\ge) \ci A\}.
  \end{equation}
\end{defn}

Note that the supremum is taken over $\ge > 0$, because $A_0 = \cj A$. The inradii $\gr_n = \gr(\tile_n)$ replace the lengths $\ell_n = 2\gr(L_n)$ of the 1-dimensional theory;
furthermore, the inradius is characterized by the following theorem.\footnote{Theorem~\ref{thm:inradius-characterization} is folkloric, but we were unable to find it in the literature and so have provided a proof.}

\begin{theorem}
  \label{thm:inradius-characterization}
  In \bRd, the inradius is the furthest distance from a point of
  $A$ to $\del A$, or the radius of the largest ball contained in
  $A$, i.e.,
  \begin{align}
    \gr(A)
      =& \sup \{\ge > 0 \suth V(A_\ge) < V(A)\}.
        \label{eqn:first inradius def}
  \end{align}
\end{theorem}
  \begin{proof}
    The continuity of the distance and volume functionals gives $V(A_\gd) < V(A)$ if and only if there is a set $U$ of positive $d$-dimensional measure contained in the interior of $A$ which is further than \gd from any point of $\del A$. For any $x \in U$, $B(x,\gd) \ci A$. Conversely, for \gd strictly less than the right-hand side of \eqref{eqn:first inradius def}, the same reasons imply the existence of the set $U$ of positive measure.
  \end{proof}

\begin{remark}
  The proof of Theorem~\ref{thm:inradius-characterization} shows that the inradius may also be defined by $\gr(A) = \sup \{dist(x,\del A) \suth \, x \in A\}$.
\end{remark}

The utility of the inradius in the present paper arises primarily
from the equality \eqref{eqn:first inradius def} and the fact
that the inradius behaves well under the action of the self-similar
system:
\begin{equation}
  \label{eqn:behavior of inradius under scaling}
  \gr_n = \gr(\tile_n) = \gr(\simt_w(\gen_q))
    = r_1^{e_1} \dots r_J^{e_J} g_q,
\end{equation}
where $r_\j$ is the scaling ratio of $\simt_\j$, and the exponent
$e_\j \in \bN$ indicates the multiplicity of the letter $\j$ in the
finite word $w \in \Wds$.

\begin{defn}
  For $q=1,\dots,Q$, the \nth[q] \emph{generating inradius} is the inradius of the \nth[q]
  generator of the tiling \tiling and denoted
  \begin{equation}
    \label{def:genir_q}
    \genir_q := \gr(\gen_q).
  \end{equation}
\end{defn}

For convenience, we may take the generators in nonincreasing order, i.e., index the generators so that
  \begin{equation}
    \label{eqn:genir ordering}
    \genir_1 \geq \genir_2 \geq \dots \geq \genir_Q >0.
  \end{equation}

\section{Measures and zeta functions} \label{ssec:Measures and zeta
functions}

In this section and the rest of the paper, any zeta function is
understood to be the 
meromorphic extension of its defining expression.

\subsection{The geometric zeta function of a fractal string}

\label{sec:The geometric zeta function of a self-similar tiling}

We recall the notion of ``fractal string'' and ``fractal spray from \cite{FGCD}, see also \cite{LaPo1,LaPo2,LaMa,La2,La3,FGNT1,HaLa}, along with \cite[Exm.~5.1 and App.~C]{La1}.

\begin{defn}\label{def:fractal-string-FGCD}
  A \emph{fractal string} is defined to be a bounded open subset of \bR, that is, a countable collection of disjoint open intervals, $L = \bigcup_{n=1}^\iy L_n$,
with lengths $\sL = \{\ell_n\}_{n=1}^\iy$. The \emph{geometric zeta
function} of such an object is
  \begin{equation}
    \label{def:string zeta fn}
    \gzL(s) = \sum_{n=1}^\iy \ell_n^s,
  \end{equation}
and can be used to study the geometry of \sL and of its (presumably
fractal) boundary $\del \sL := \del L$.
\end{defn}
Observe that $\gzL(s)$ is the Mellin transform of the measure
  \begin{equation}
    \label{def:ghL}
    \ghL = \sum_{n=1}^\iy \gd_{1/\ell_n},
  \end{equation}
where $\gd_x$ denotes the Dirac mass (or Dirac measure) at $x$.
Thus,
  \begin{equation}
    \label{def:zeta fn as Mellin transf of genl string}
    \gzL(s) = \int_0^\iy x^s \, d\ghL(x).
  \end{equation}

The following definition first appeared in \cite{LaPo2}.

\begin{defn}\label{def:fractal-spray}
  Let $\gen \ci \bRd$ be a nonempty bounded open set, which we will call the \emph{generator} (or \emph{basic shape}). Then a \emph{fractal spray} \tiling is a bounded open subset of \bRd which is the disjoint union of open sets $\gW_n$ for $n=1,2,\dots$, where each $\gW_n$ is congruent to $\ell_n \gen$, the homothetic of \gen by $\ell_n$. Here, $\sL = \{\ell_n\}_{n=1}^\iy$ is a fractal string and we say that \tiling is ``scaled by'' \sL.
  It is clear that a self-similar tiling as discussed in Definition~\ref{def:ss tiling} is a special case of fractal spray.
\end{defn}
Thus, any fractal string can be thought of as a fractal spray on the generator $\gen =(0,1)$, the unit interval. In the context of the current paper, a self-similar tiling is a union of fractal sprays on the generators $\gen_1,\dots,\gen_Q$, each scaled by a fixed \emph{self-similar} string. A general fractal spray may have multiple generators, as long as they are all scaled by the same 
fractal string \sL. However, for the remainder of this paper we consider only a single generator \gen. Indeed, as mentioned in \S\ref{ssec:Tilings with multiple generators}, the multiple-generator case can readily be reduced to the case of a single generator.

\begin{remark}\label{rem:generalized-strings}
  One can also define a \emph{generalized fractal string} \gh to be a locally finite regular Borel measure on $\bR_+$ with support bounded away from 0, as is done in \cite[\S4]{FGCD}. With regard to scaling, this is perhaps the greatest degree of generality in which the results of the present paper also hold. Indeed, all the main results of this paper remain true when \ghs is replaced by a generalized fractal string \gh, including the technical content of the appendices. However, since the geometric realization of a generalized fractal string (or spray) is unclear at this point, we describe tube formulas and other concepts from geometric measure theory in terms of \ghs.
\end{remark}

\subsection{The scaling zeta function}
\label{ssec:scaling-zet-fn}

We now extend Definition~\ref{def:fractal-string-FGCD} to higher dimensions. In 1 dimension, the length of an interval is just twice its inradius, and the distinction between the scale of a set and its volume is blurred. In higher dimensions, the two are related by a power law, i.e., in \bRd, the volume of a compact set of full dimension will change by a factor of $r^d$ when the set is scaled by a factor $r$. In \cite{FGCD}, the geometric zeta function \gzL records the measures of all the open sets comprising the string. However, in higher dimensions it is easier to use a zeta function (\gzs, introduced just below) to record the scales (inradii) of the open sets in question.  

In Definition~\ref{def:geometric zeta fn of a fractal spray}, we will introduce the \tubezf zeta function of the tiling \gzT, which is the higher-dimensional analogue of \gzL, and is defined in terms of \gzs and the generators. The tiling zeta function \gzT encodes the density of geometric states of \tiling and acts as a generating function for the geometry of the entire tiling. The scaling zeta function encodes only scaling data (inradii). In \cite{FGCD}, both of these roles are essentially played by \gzL. 

\begin{defn}
  \label{defn:scaling-and-geom-meas}
  For a fractal spray scaled by a given fractal string $\sL = \{\ell_n\}_{n=1}^\iy$, we use the inradii $\gr_n := \frac{\ell_n}{2}$ to define the \emph{scaling measure} by 
  \begin{align}
    \label{def:scaling-measure-spray}
    \ghs(x) := \sum_{n=1}^\iy \gd_{1/\gr_n}(x).
  \end{align}
\end{defn}

\begin{defn}
  \label{defn:scaling and geom zetas}
  The \emph{scaling zeta function} is the Mellin transform of the scaling measure:
  \begin{align}
    \label{def:scaling-zeta}
    \gzs(s)
      := \int_0^\iy x^{-s} \, d\ghs(x).
  \end{align}
\end{defn}

\begin{remark}\label{rem:self-similar-scaling-zeta}
  In the special case of a self-similar tiling, as in \S\ref{sec:The self-similar tiling}, it is easy to see from Definition~\ref{def:fractal-spray} that the scaling measure encodes the scaling factors of the self-similar system \simt as a sum of Dirac masses:
  \begin{align}
    \label{def:scaling-measure-tiling}
    \ghs(x) := \sum_{w \in \Wds} \gd_{1/r_w}(x).
  \end{align}
  In this case, one can see that the scaling zeta function \gzs encodes the combinatorics of the scaling ratios $\{r_j\}_{j=1}^J$ of $\simt=\{\simt\}_{j=1}^J$, and is thus a generating function for the scaling properties of \simt:
  \begin{align}
   \label{eqn: self-similar-scaling-zeta-function}
    \gzs(s) 
       = \sum_{w \in \Wds} r_w^s
       = \sum_{k=0}^\iy \sum_{w \in \Wds_k} r_w^s.     
    \end{align}
    However, there is another, simpler (and more useful) form of \gzs, which we now explain. 
\end{remark}

Theorem~\ref{thm:scaling-zeta-fn-simplified} is the higher-dimensional counterpart of \cite[Thm.~2.4]{FGCD}, and can, in fact, be viewed as a corollary of it; see \S\ref{ssec:comparison to FGNT}. Indeed, it is proved in precisely the same way.

\begin{theorem}
  \label{thm:scaling-zeta-fn-simplified}
  The scaling zeta function of a self-similar tiling is
  \begin{align}
    \label{eqn:scaling zeta function simplified}
    \gzs(s) 
       = \frac1{1-\sum_{\j=1}^J r_\j^s}.
  \end{align}
  This remains valid for the meromorphic extension of \gzs to all of \bC.
\end{theorem}

\begin{defn}
  \label{defn:complex dimensions}
  We can now define the \emph{scaling (complex) dimensions} of a tiling \tiling as the poles of the scaling zeta function:
  \begin{align}
    \label{eqn:complex dimensions def}
    \Ds := \{\gw \in \bC \suth \gzs(s) \text{ has a pole at }
    \gw\}.
  \end{align}
\end{defn}

\subsection{Comparison with \cite{FGCD}}
  \label{ssec:comparison to FGNT}

  Although the measure and zeta function introduced in Definition~\ref{defn:scaling-and-geom-meas} and Definition~\ref{defn:scaling and geom zetas} above correspond to fractal subsets of \bRd, it is crucial to note that (setting geometric interpretation aside) they are formally identical to the objects \sL (thought of as a measure) and \gzL studied in \cite{FGCD}. To be precise, the scaling measure \ghs is a fractal string of the sort studied in Chap.~2--3 of \cite{FGCD}, and a generalized fractal string of the kind introduced in Chap.~4 of \cite{FGCD}. Consequently, all of the explicit formulas developed in \cite{FGCD} are applicable to the measures and zeta functions described in the present paper. This is key to the proof of Theorem~\ref{thm:fractal-spray-tube-formula}, the tube formula for fractal sprays.
  
  In the special case when \ghs is the scaling measure of a self-similar tiling, \gzs may be thought of as the geometric zeta function of a self-similar string (as in \cite{FGCD}) with scaling ratios $\{r_j\}_{j=1}^J$ and a single \emph{gap}\footnote{In this paper, we use the term ``generator'' in place of ``gap''.}, which has been normalized so as to have $\ell_1=1$, where $\ell_1$ is the first length in the string. 
  Let $D$ be the unique real number satisfying $\sum_{j=1}^J r_j^D = 1$. One can check (as in \cite{FGCD}, \S5.1) that for some real constant $c > D$,
  \begin{align}
    \label{eqn:measure-zeta reciprocity}
    \ghs(x) &= \frac1{2\gp \ii} \int_{c-\ii\iy}^{c+\ii\iy}  x^{s-1} \gzs(s) \, ds,
    \text{\; \, and \;}
    \gzs(s) = \int_0^\iy x^{-s} \, \ghs(dx).
  \end{align}

  Additionally, the structure theorem for complex dimensions of self-similar strings 
  \cite[Thm.~3.6]{FGCD} (in the special case of a single gap), holds for the set of scaling complex dimensions of Definition~\ref{defn:complex dimensions}.
  By \eqref{eqn:scaling zeta function simplified}, \Ds consists of the set of complex solutions of the complexified Moran equation $\sum_{j=1}^J r_j^s = 1$ which is studied in detail in \cite[Chap.~2--3]{FGCD}. In particular, the scaling complex dimensions lie in a horizontally bounded strip of the form $D_\ell \leq \Re s \leq D$, where $D$ is as just above and $D_\ell < D$ is some other (finite, possibly negative) constant. The positive number $D$ is called the \emph{similarity dimension} of \simt (or of its attractor \attr) and coincides with the abscissa of convergence of \gzs \cite{FGCD}, Thm. 1.10.\footnote{If the self-similar system defining \attr satisfies the `open set condition' (see \cite{Hut}, as described in \cite{Fal1} or \cite{Kig}), as when the tileset condition is satisfied, then $D$ coincides with the Hausdorff and Minkowski dimensions of \attr.} Furthermore, the following dichotomy prevails:
  \begin{itemize}
    \item \emph{Lattice case}. When the logarithms of the scaling ratios $r_j$ are each an integer power of some common positive real number, the scaling complex dimensions lie periodically on finitely many vertical lines, including the line $\Re s = D$. In this case, there are infinitely many complex dimensions with real part $D$.

    \item \emph{Nonlattice case}. Otherwise, the scaling complex dimensions are quasiperiodically distributed and $s=D$ is the only complex dimension with real part $D$. However, there exists an infinite sequence of scaling complex dimensions approaching the line $\Re s=D$ from the left. In this case, the set $\{\Re s \suth s \in \D\}$ appears to be dense in (finitely many compact subintervals of) $[D_\ell,D]$.
  \end{itemize}
  It has been proven in \cite{FGCD} that for $d=1$, the
  attractor of \simt fails to be Minkowski measurable if and only if
  \gzL has nonreal complex dimensions with real part $D$, and in
  \cite{FGCD}, Conj.~12.18, this is conjectured to hold also in
  higher dimensions. See also
  Remark~\ref{rem:dimensions-and-measurability}.


\section{The Generators}
\label{sec:tube formula for tilings}

The inner tube formula for a self-similar tiling will consist of the sum
of the inner tube formulas for each tile, and each of these can
be expressed as a rescaled version of the tube formula for a
generator. That is, if $\tile = \simt_w(G)$ for some $w \in
\Wds$, then the inradius of such a tile is
\[\gr = \gr(\tile) = \gr(\simt_w(\gen)) = r_w \genir = r_1^{e_1} \dots r_J^{e_J} \genir,\]
and invariance of Lebesgue measure under rigid motions gives
\begin{align}
  \label{eqn:V for a tile}
  V(\tile,\ge)
  &= V(\simt_w(\gen),\ge)
   = V(r_w \gen,\ge).
\end{align}
Thus, it behooves us to find an expression for
\begin{align}
  \gtf(x,\ge) := V(\tfrac1x\gen,\ge),
    \label{def:gtf}
\end{align}
where $\frac1x \gen$ is a homothetic image of \gen, scaled by
some factor $\frac1x>0$. Then $\gtf(1,\ge)$ gives the inner
tube formula for $\gen$, and $\gtf(x,\ge)$ is the volume of a
tile which is similar to \gen but which has been scaled by
$1/x$. The motivation for defining \gtf in terms of $1/x$
(rather than $x$) appears in \eqref{eqn:split integral}. 
  Since $\gtf(x,\ge)$ gives the inner tube formula for a set congruent to  $(1/x)\gen$, it is useful for computing the inner tube formula for a general fractal spray.

\begin{defn}
  \label{def:Steiner-like}
  A \emph{Steiner-like representation}%
  \footnote{We are grateful to Steffen Winter for pointing out that \eqref{eqn:def-prelim-Steiner-like-formula} is really a property of (a choice of) the tube formula, as opposed to a property of the generator, and for suggesting the term ``Steiner-like representation'' (to replace the term ``Steiner-like generator'', which appeared in a previous draft).} 
  of a generator \gen is an inner tube formula  of the form
  \begin{align}\label{eqn:def-prelim-Steiner-like-formula}
    V(\gen,\ge) = \sum_{k=0}^{d} \crv_k(\gen,\ge) \ge^{d-k}, 
    \qq \text{for } \ge < \genir, 
  \end{align}
 where each $\crv_k(\gen,\ge)$ is some
  reasonably nice (e.g., bounded and locally integrable) function for $\ge \in [0, \genir)$.
  In particular, we require that
  \begin{enumerate}[(i)]
    \item each $\crv_k(\gen,\ge)$ is homogeneous of degree $k$, so that for
    $\gl > 0,$
      \label{itm:homogeneity of crv}
      \begin{align}
        \label{eqn:homogeneity of crv}
        \crv_k\left(\gl\gen,\gl\ge\right) = \crv_k(\gen,\ge) \, \gl^k,
        \q\text{and}
      \end{align}
    \item each $\crv_k(\gen,\ge)$ is rigid motion invariant, so that
      \label{itm:translation invariance of crv}
      \begin{align}
        \label{eqn:translation invariance of crv}
        \crv_k\left(T(\gen),\ge\right) = \crv_k(\gen,\ge),
      \end{align}
      for any (affine) isometry $T$ of \bRd.
    \item for each $\crv_k(\gen,\ge)$, $k=0,1,\dots,d$, the
      limit $\lim_{\ge \to 0^+} \crv_k(\gen,\ge)$ exists in \bR.
      \label{itm:integrability of crv}
  \end{enumerate}
\end{defn}
We have chosen the term ``Steiner-like'' for Definition~\ref{def:Steiner-like} because the intrinsic volumes $\gm_k$ (see Theorem~\ref{thm:Steiners-formula-basic} and the ensuing discussion) satisfy the following properties:
\begin{enumerate}[(i)]
  \item each $\gm_k$ is homogeneous of degree $k$, so that for $x >0,$
  \label{itm:homogeneity of mu}
  \begin{align}
    \label{eqn:homogeneity of mu}
    \gm_k\left(x A\right) = \gm_k(A) \, x^k, \text{ and}
  \end{align}
\item each $\gm_k(A)$ is rigid motion invariant, so that
  \label{itm:translation invariance of mu}
  \begin{align}
    \label{eqn:translation invariance of mu}
    \gm_k\left(T(A)\right) = \gm_k(A),
  \end{align}
  for any (affine) isometry $T$ of \bRd.
\end{enumerate}
The description of $\crv_k$ given in the conditions of Definition~\ref{def:Steiner-like} is intended to emphasize the resemblance between $\crv_k$ and $\gm_k$. However, $\crv_k$ \emph{may be signed} (even when \gen is convex and $k=d-1$ or $d$) and is more complicated than $\gm_k$ in general. In contrast, the Federer curvature measures $\gQ_k$ are always positive for convex sets; cf.~\cite{Fed} and \cite{Schn2}.

\emph{Note added in proof:} Property (iii) of Definition~\ref{def:Steiner-like} eventually proved to be unnecessary, and will be omitted in future work, e.g., \cite{Pointwise}.

\subsection{The tube formula for a monophase generator}
\label{ssec:The tube formula for a generator}

In this paper, we treat only the special case of a generator \gen which has a Steiner-like representation whose coefficient functions $\crv_k(\gen,\ge)$ are constant functions of \ge for $0 \leq \ge < \genir$, in which case \gen is called \emph{monophase} as in Definition~\ref{def:monophase}. 
We treat more general generators in the forthcoming collaboration with Steffen Winter \cite{Pointwise}. 

\begin{defn}
  \label{def:monophase}
  A generator \gen is said to be a \emph{monophase generator}, or to have a \emph{monophase tube formula}, iff one can write
  \begin{align}
    V(\gen,\ge) = \gtf(1,\ge) = \sum_{k=0}^{d-1} \crv_k(\gen) \ge^{d-k},
    \qq\text{for $\ge < \genir$,}
      \label{eqn:def-prelim-monophase-formula}
  \end{align}
  for some constants $\crv_k(\gen) \in \bR$, $k=0,1,\dots,d-1$.
\end{defn}
Not every polyhedral generator \gen is monophase; the more general pluriphase case is discussed in \S\ref{ssec:More general generators}. In general, the computation of $V(\gen,\ge)$ may be nontrivial. 
So far, we have only defined $V(\gen,\ge)$ for $\ge < \genir$. To extend
  it to all of $\bR^+$, note that $V(\gen,\ge)$ is just the
  Lebesgue measure of \gen for $\ge \geq \genir$. Therefore,
  define
  \begin{align*}
    \crv_k(\gen,\ge) &:= \crv_k(\gen) \gc_{[0,\genir)}(\ge),
      &&i=0,1,\dots,d-1 \notag \\
    \crv_d(\gen,\ge)  &:= -\crv_d(\gen) \gc_{[\genir,\iy)}(\ge),
  \end{align*}
  where $\crv_k(\gen)$ is as in
  \eqref{eqn:def-prelim-monophase-formula}, $\gc_A$ is the usual
  characteristic function of the set $A$, and $\gk_d(\gen)$ is defined to be the negative of the $d$-dimensional Lebesgue measure of \gen:
\begin{equation}\label{eqn:def:kappad}
  \crv_d(\gen) = -\vol[d](\gen).
\end{equation}
Now we have
  \begin{align}
    V(\gen,\ge) = \sum_{k=0}^{d} \crv_k(\gen,\ge) \ge^{d-k},
    \qq\text{for $\ge \geq 0$.}
      \label{eqn:def-monophase-formula}
  \end{align}

\begin{remark}\label{rem:no-d-in-monophase}
  Note that for the case $k=d$, one must have $\lim_{\ge \to 0^+} \crv_d(\gen,\ge) = 0$. This implies that $V(\gen,\ge)$ can have no constant term when \gen is monophase, and hence the sum in \eqref{eqn:def-prelim-monophase-formula} does not include a \nth[d] term. (The same may not hold when \gen is pluriphase.)
\end{remark}

\begin{theorem}\label{thm:monophase-formula}
  If \gen is monophase, then for any tile congruent to the
  homothetic image $\frac1x \gen$, the inner tube formula is given by
  \begin{equation}
    \label{eqn:thm:monophase-formula}
    \gtf(x,\ge) =
    \begin{cases}
      \sum_{k=0}^{d-1} \crv_k(\gen) x^{-k} \ge^{d-k}, &\ge \leq \genir/x, \\
      -\crv_d(\gen)x^{-d}, &\ge \geq \genir/x.
    \end{cases}
  \end{equation}
\end{theorem}
\begin{proof}
  From \eqref{eqn:def-monophase-formula}, we have $V(\gen,\ge) = \gtf(1,\ge)$; we would like to adapt this formula so as to obtain a
  tube formula valid for a tile of any size. Note that
  $V_{r\gen}(r\ge) = r^d V(\gen,\ge)$, as both expressions are
  measuring congruent regions in \bRd. Hence for $\ge < \genir$, one has
  \begin{align*}
    \sum_{k=0}^{d-1} r^k \crv_k(\gen, \ge) (r\ge)^{d-k}
    &= r^d V(\gen,\ge)
    = V({r\gen},r\ge)
    = \sum_{k=0}^{d-1} \crv_k(r\gen, r\ge) (r\ge)^{d-k},
  \end{align*}
  and thus for $\ge < \genir/x$, one has
  \begin{align*}
    \gtf(x,\ge) = V({\tfrac1x \gen},\ge)
    =  \sum_{k=0}^{d-1} \left(\tfrac1x\right)^k \crv_k(\gen, x\ge) \ge^{d-k}
    =  \sum_{k=0}^{d-1} \crv_k(\gen, x\ge) x^{-k} \ge^{d-k}.
  \end{align*}
  Since $\crv_k(\gen; x\ge) = \crv_k(\gen) \gc_{[0,\genir)}
  (x\ge) = \crv_k(\gen) \gc_{[0,\genir/x)} (\ge)$ for $k = 0,1,\dots,d-1$,
  it is clear that \eqref{def:gtf} may be expressed as
  \begin{equation}
    \label{eqn:monophase-form-of-gtf}
    \gtf(x,\ge) = \sum_{k=0}^{d} \crv_k(\tfrac1x\gen, \ge) \ge^{d-k} =
    \begin{cases}
      \sum_{k=0}^{d-1} \crv_k(\gen) x^{-k} \ge^{d-k}, &\ge \leq \genir/x, \\
      -\crv_d(\gen)x^{-d}, &\ge \geq \genir/x,
    \end{cases}
  \end{equation}
  where the constants $\crv_k(\gen)$ are as defined in
  \eqref{eqn:def-prelim-monophase-formula} for $k = 0,1,\dots,d-1$,
  and in \eqref{eqn:def:kappad} for $k=d$.
\end{proof}

The function $\gtf(x,\ge)$ gives the volume of the \ge-neighbourhood of
a tile which is congruent to a generator scaled by $1/x$. The
value $\ge=\genir/x$ corresponds to the value of \ge at which
the inner \ge-neighbourhood of the tile becomes equal to the tile
itself. Thus, $\gtf$ is continuous (but generally not differentiable) at
$\ge=\genir/x$.

\subsection{More general generators}
\label{ssec:More general generators}

Not every generator \gen is monophase, so we introduce the pluriphase case in Definition~\ref{def:pluriphase}. The most general case is discussed in \cite{Pointwise}. In fact, even if \gen is polyhedral or convex, it still may not be monophase. Example~\ref{exm:pluriphase_generator} gives a set which is convex and pluriphase but not monophase and the Cantor Carpet discussed in \cite[Ex.~4.1]{Pointwise} gives a tiling with a nonconvex generator (in the shape of a Swiss cross) which is polyhedral but not even pluriphase.


\begin{defn}
  \label{def:pluriphase}
  A generator \gen is said to be a \emph{pluriphase} generator iff its inner tube formula is given by a piecewise polynomial function of \ge. Equivalently, \gen is pluriphase iff the functions $\crv_k(\gen;\ge)$ of \eqref{eqn:def-prelim-Steiner-like-formula} are piecewise constant for $k=0,1,\dots,d$. 
\end{defn}
  
It is possible (though doubtful) that all convex generators are pluriphase, but this has not yet been proved. However, it seems likely that all convex polyhedra are pluriphase. For situations even more general, it is an interior version of Federer's notion of \emph{reach} (see \cite{Fed}) that is required. For such cases, the inner tube formula will be obtained in \cite{FCM} via the more general methods of \cite{HuLaWe} and others. 
just like Theorem~\ref{thm:monophase-formula}.

\begin{exm}[A pluriphase generator]
  \label{exm:pluriphase_generator}
  Consider a fractal spray on a generator \gen consisting of a
  $2\times2$ square with one corner replaced by a circular arc,
  as depicted in Figure~\ref{fig:pluriphase}. This generator has
  inradius $\genir = \gr(\gen) = 1$ and is pluriphase, but not
  monophase. Indeed, the relevant partition is
  \begin{equation}\label{eqn:exm:pluriphase-partition}
    \{0=\ge_0, \ge_1=1/2, \ge_2=1\},
  \end{equation}
  and the tube formula for \gen is
  \begin{equation}\label{eqn:exm:pluriphase-gtf}
    \gtf(1,\ge) =
    \begin{cases}
      (8+\tfrac\gp4)\ge - (5+\tfrac\gp4)\ge^2, &\ge_0 \leq \ge \leq \ge_1 \\
      \frac{\gp-4}{16} + 8\ge - 4\ge^2, &\ge_1 \leq \ge \leq \ge_2 \\
      \frac{\gp-4}{16} + 4, &\ge_2 \leq \ge .
    \end{cases}
  \end{equation}

  \begin{figure}
    \centering
    \scalebox{0.70}{\includegraphics{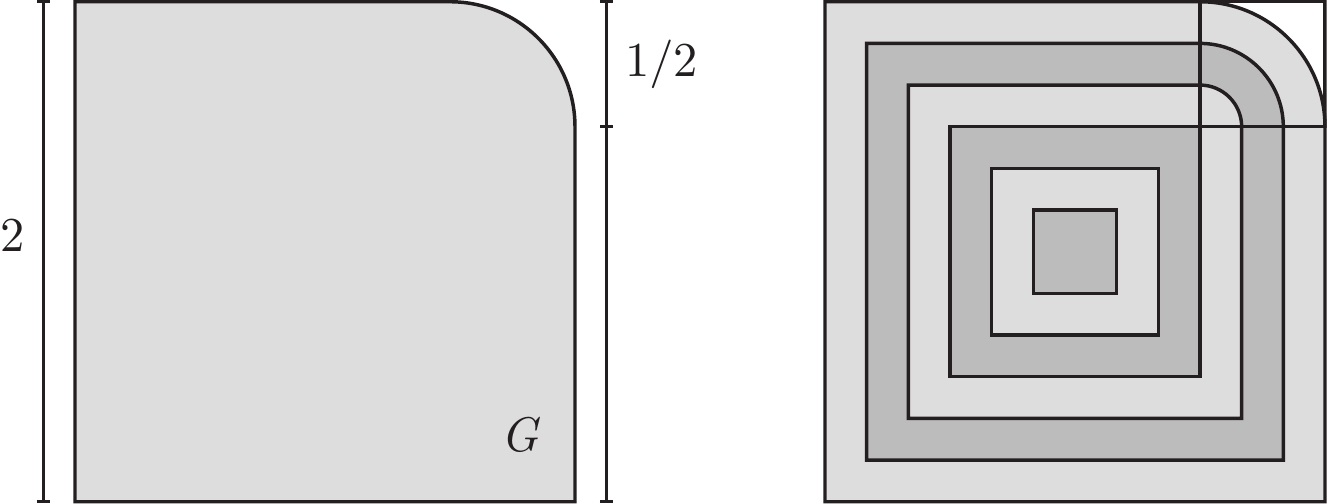}}
    \caption{A pluriphase generator which is not monophase.}
    \label{fig:pluriphase}
  \end{figure}
\end{exm}

\subsection{Tilings with one generator}
\label{ssec:Tilings with one generator}

Suppose we have a tiling \tiling with just one generator \gen.
Then the inner tube formula of $\tiling$ is given by
\begin{align}
  V(\tiling,\ge)
  &= \sum_{n=1}^\iy V({\tile_n},\ge)
  = \sum_{\gr_n \geq \ge} V({\tile_n},\ge) + \sum_{\gr_n < \ge} V({\tile_n},\ge),
  \label{eqn:split sum}
\end{align}
much as in \cite{LaPo2}, Eqn.~(3.2). Recall that $\gr_n$ is the inradius of the tile $\tile_n$. For $\tile_n = \simt_w(\gen_q)$, invariance under rigid motions allows us to use the equality \eqref{eqn:V for a tile} to rewrite the sums in \eqref{eqn:split sum} as integrals with respect to the scaling measure \ghs:
\begin{align}
  V(\tiling,\ge)
  &= \sum_{\gr_n^{-1} \leq 1/\ge} V({\tile_n},\ge) + \sum_{\gr_n^{-1} > 1/\ge} V({\tile_n},\ge) \notag \\
  &= \int_0^{\genir/\ge} V({\tfrac1x\gen},\ge) \,d\ghs(x)
    + \gm_d(\gen) \int_{\genir/\ge}^\iy x^{-d} \,d\ghs(x)
    \label{eqn:split integral} \\
  &= \int_0^\iy \gtf(x,\ge) \, d\ghs(x) \notag \\
  &= \la \ghs, \gtf \ra, \vstr
  \label{def:action of a measure on a test fn}
\end{align}
where \gtf is a `test function' giving the volume of a tile which is similar to \gen, but which has been scaled by a factor of $1/x$; see \eqref{def:gtf}. Although \gtf is not smooth, it fits the criteria given in Theorem~\ref{thm:ext dist explicit formula} and is thus amenable to the distributional techniques of \cite[\S5.4]{FGCD}.

\subsection{Tilings with multiple generators}
\label{ssec:Tilings with multiple generators}

Upon replacing $\gen$ by $\gen_q$, we use the notation $V_q, \gtfq, \crv_{qk}$, etc., to refer to the corresponding quantity for the \nth[q] generator. For example, $\gtf(x,\ge)$ is replaced by $\gtfq(x,\ge) = \gtf[\gen_q](x,\ge)$, the volume of the \ge-neighbourhood of a tile which is similar to $\gen_q$ but which has been scaled by $x$.

The contribution to $V(\tiling,\ge)$ resulting from one generator $\gen_q$ and its successive images is $V(\gen_q,\ge) := \la \ghs, \gtfq \ra$, so the case of multiple generators can be reduced to a sum of single-generator tilings via the formula
\begin{equation}
  \label{eqn:V as eta acting on v_e}
  V(\tiling,\ge)
  = \sum_{q=1}^Q V(\gen_q,\ge)
  = \sum_{q=1}^Q \la \ghs, \gtfq \ra.
\end{equation}
For a concrete example of how this is done, see the example of the pentagasket in \S\ref{exm:Pentagasket Tiling}.

Henceforth, we will always assume there is only a single generator, as this simplifying assumption will clarify the exposition. The single generator will always be denoted by \gen.


\section{Distributional Explicit Formulas for
Fractal Strings} \label{sec:The distributional theory of fractal
strings}


These four definitions and the three theorems that follow them are adapted from
\cite[\S5.3]{FGCD}. The technical details described here are used in the proof of Theorem~\ref{thm:fractal-spray-tube-formula}, especially in Appendix~\ref{app:Rigourous Definition of gzT} and Appendix~\ref{app:The Error Term and Estimate}. The reader can easily skim or skip this section on a first reading.



\begin{defn}
  \label{def:screen}
  Let $f:\bR \to \bR$ be a bounded Lipschitz continuous function. Then the \emph{screen} is $S = \{f(t) + \ii t \suth t \in \bR\}$, the graph of a function with the axes interchanged. We let
  \begin{align}
    \label{eqn:infS}
    \inf S &:= \inf\nolimits_t f(t) = \inf\{\Re s \suth s \in S\}, \text{ and} \\
    \label{eqn:supS}
    \sup S &:= \sup\nolimits_t f(t) = \sup\{\Re s \suth s \in S\}.
  \end{align}
  The screen is thus a vertical contour in \bC. The region to the right of the screen is the set $W$, called the \emph{window}:
  \begin{align}
    \label{eqn:window}
    W &:= \{z \in \bC \suth \Re z \geq f(\Im z)\}.
  \end{align}
  The poles of \gzs  are called the \emph{scaling dimensions} and those which lie in the window are called the \emph{visible scaling dimensions}; the set of them is denoted
  \begin{align}
    \label{eqn:def:visible-scaling-dimns}
    \Ds(W) = \{\gw \in W \suth \lim_{s \to \gw} |\gzs (s)| = \iy\}.
  \end{align}
\end{defn}

\begin{defn}
  \label{def:languid}
  The 
  scaling measure \ghs (as in Definition~\ref{defn:scaling-and-geom-meas})
  is said to be \emph{languid} if its associated zeta function \gzs  satisfies certain growth conditions relative to the screen. 
  Specifically, let $\{T_n\}_{n \in \bZ}$ be a sequence in \bR such that $T_{-n} < 0 < T_n$ for $n \geq 1$, and
  \begin{align}
    \label{eqn:Tn seq conds}
    \lim_{n \to \iy} T_n = \iy,
    \lim_{n \to \iy} T_{-n} = -\iy, \text{ and }
    \lim_{n \to \iy} \frac{T_n}{|T_{-n}|} = 1.
  \end{align}
  For \ghs to be languid, there must exist real constants $\languidOrder,c>0$ and a sequence $\{T_n\}$ as described in \eqref{eqn:Tn seq conds}, such that\\
  \vstr \textbf{L1} \; For all $n \in \bZ$ and all $\gs \geq f(T_n)$,
  \begin{align}
    \label{eqn:L1}
    |\gzs (\gs + \ii T_n)| \leq c \cdot \left(|T_n| + 1\right)^\languidOrder,
    \text{ and}
  \end{align}
  \textbf{L2} \; For all $t \in \bR$, $|t| \geq 1$,
  \begin{align}
    \label{eqn:L2}
    |\gzs (f(t) + \ii t)| \leq c \cdot |t|^\languidOrder.
  \end{align}
  In this case, \ghs is said to be \emph{languid of order \languidOrder}.
\end{defn}

\begin{defn}
  \label{def:strongly languid}
  The 
  scaling measure \ghs is said to be \emph{strongly languid} if it satisfies \textbf{L1} and the condition \textbf{L2'}, which is clearly stronger than \textbf{L2}:\\
  \vstr
  \textbf{L2'} \; There exists a sequence of screens $S_m(t) = f_m(t)+ \ii t$ for $m \geq 1$, $t \in \bR$, with $\sup S_m \to -\iy$ as $m \to \iy$, and with a uniform Lipschitz bound. Additionally, there must exist constants $A,c>0$ such that
  \begin{align}
    \label{eqn:L2'}
    |\gzs (f(t) + \ii t)| \leq c \cdot A^{|f_m(t)|}(|t|+1)^\languidOrder,
  \end{align}
  for all $t \in \bR$ and $m \geq 1$.
\end{defn}

It turns out that Definition~\ref{def:strongly languid} is always satisfied when \ghs is the scaling measure of a self-similar tiling; see \cite[\S6.4]{FGCD}.

Taking \cite[Thm.~5.26 and Thm.~5.30]{FGCD} at level $k=0$ gives the following distributional explicit formula for the action of a scaling measure \ghs on a test function $\testfn \in C^\iy(0,\iy)$.
Note that \testfn may not have compact support; only the decay properties \eqref{eqn:test fn first req for extended formula}--\eqref{eqn:test fn second req for extended formula} are required.

\begin{theorem}[Extended distributional explicit formula]
  \label{thm:ext dist explicit formula}
  Let \ghs be a 
  scaling measure which is \slow of order \languidOrder. Let $\testfn \in C^\iy(0,\iy)$ with \nth derivative satisfying, for some $\gd > 0$, and every integer $n \in \{0,1,\dots,N = [\languidOrder]+2\}$,
  \begin{equation}
    \label{eqn:test fn first req for extended formula}
    \testfn^{(n)}(x) = O\left(x^{-n-D-\gd}\right) \q \text{ as } x \to \iy,
    \q \text{and}
  \end{equation}
  \begin{equation}
    \label{eqn:test fn second req for extended formula}
    \testfn^{(n)}(x) = \sum_\ga a_\ga^{(n)} x^{-\ga-n} + O\left(x^{-n-\inf S + \gd}\right)
    \q \text{ as } x \to 0^+.
  \end{equation}
  Then we have the following distributional explicit formula for \ghs:
  \begin{equation}
    \label{eqn:action of dist on test fn}
    \left\la \ghs,\testfn\right\ra
    = \sum_{\gw \in \Ds} \negsp[4] \res{\gzs (s)\testfnt(s)}
     + \negsp[10]\sum_{\ga \in W\less\Ds} \negsp[15] a_\ga^{(0)} \gzs (\ga)
     + \left\la \R, \testfn\right\ra,
  \end{equation}
  where the error term $\R(x)$ is the distribution given by
  \begin{equation}
    \label{def:ext-dist-error-term}
    \left\la \R, \testfn\right\ra
    = \tfrac1{2\gp \ii} \int_S \gzs (s) \testfnt(s)\,ds
  \end{equation}
  and estimated by
  \begin{equation}
    \label{def:error term estimate}
    \R(x) = O\left(x^{\sup S - 1}\right), \q \text{ as } x \to \iy.
  \end{equation}
\end{theorem}

  Here, \testfnt is the Mellin transform of the function \testfn, defined by
  \begin{equation}
    \label{def:Mellin transform}
    \testfnt(s) := \int_0^\iy x^{s-1} \testfn(x) \,dx.
  \end{equation}

Note: the sum in (\ref{eqn:test fn second req for extended
formula}) is over finitely many complex exponents \ga with $\Re
\ga > -\gs_l+\gd$; we express this by saying that
\testfn has an asymptotic expansion of order $-\gs_l + \gd$ at
0.

Taking \cite{FGCD}, Thm.~5.27, at level $k=0$ gives the following distributional explicit formula for the action of a scaling measure \ghs on a test function \gy. Note that in addition to requiring $\testfn \in C^\iy(0,\iy)$, we now also require that \gy is a finite linear combination of terms $x^{-\gb}e^{-c_\gb x}$ in a neighbourhood of the interval $(0,A]$, where $A$ is the same constant as in Definition~\ref{def:strongly languid}.

\begin{theorem}[Extended distributional formula, without error term]
\label{thm:ext_distributional_without error} Let~$\ghs$ be a strongly languid 
scaling measure. Let $q \in \bN$ be such that $q>\max\{ 1,\languidOrder\},$ where \languidOrder is as in Definition~\ref{def:languid}. Further, let \testfn be a test function that is $q$ times continuously differentiable on $(0,\infty )$. Assume that the \nth[j] derivative $\testfn^{(j)}(x)$ satisfies \eqref{eqn:test fn first req for extended formula} and \eqref{eqn:test fn second req for extended formula}, and that there exists a $\gd > 0$ such that
\begin{align}
  \testfn^{(j)}(x)
  &= \sum_\ga a^{(j)}_\ga x^{-\ga} e^{-c_\ga x},
  \qq \text{for } x \in (0,A+\gd), 0 \leq j \leq q.
\end{align}
Then formula \eqref{eqn:action of dist on test fn} holds with $\sR
\equiv 0$.
\end{theorem}

\begin{theorem}[Tube formula for fractal strings {\cite{FGCD}, Thm.~8.1}]
  \label{thm:tube formula for fractal strings}
  Let $\sL = \{\ell_n\}_{n=1}^\iy$ be a fractal string with 
  \slow zeta function $\gzL = \sum_{n=1}^\iy \ell_n^s$. 
  Then the volume of the (one-sided) tubular neighbourhood of radius \ge of the boundary of \sL is given by the following distributional explicit formula for test functions $\testfn \in C_c^\iy(0,\iy)$, the space of $C^\iy$ functions with compact support contained in $(0,\iy)$:
  \begin{equation}
    \label{eqn:1-dim tube formula}
    V(\sL,\ge)
    = \sum_{\gw \in \D_\sL(W)} \res{\frac{\gzL (s)(2\ge)^{1-s}}{s(1-s)}}
     + \{2\ge \gzL (0)\} + \sR(\ge).
  \end{equation}
  Here the term in braces is only included if $0 \in W \less \Ds(W)$, and
  $\sR(\ge)$ is the error term, given by
  \begin{equation}
    \label{eqn:1-dim error term}
    \sR(\ge)
    = \frac1{2\gp \ii} \int_S \frac{\gzL (s) (2\ge)^{1-s}}{s(1-s)} \, ds
  \end{equation}
  and estimated by
  \begin{equation}
    \label{eqn:1-dim error bound}
    \sR(\ge) = O(\ge^{1-\sup S}), \qq \text{as } \ge \to 0^+.
  \end{equation}
\end{theorem}
The meaning of \eqref{def:error term estimate} and \eqref{eqn:1-dim
error bound}, the order of the distributional error term, is given in
Definition~\ref{def:order of a dist} of Appendix~\ref{app:The Error Term
and Estimate}.
%
When \sL is a \emph{self-similar} fractal string, the results of Theorem~\ref{thm:tube formula for fractal strings} may be strengthened as described in \S\ref{ssec:The self-similar case} and in \cite[\S8.4]{FGCD}. In particular, one may take $W= \bC$ and $\sR \equiv 0$.


\section{The Tube Formula for Fractal Sprays}
\label{sec:The tube formula for self-similar tilings}

We now present the main result of the paper, a higher-dimensional analogue of Theorem~\ref{thm:tube formula for fractal strings}. While the proof is similar in spirit to the proof of the tube formula for fractal strings obtained in \cite[\S8.1]{FGCD} (cf. Theorem~\ref{thm:tube formula for fractal strings}), it is significantly more involved, especially if Appendices \ref{app:Rigourous Definition of gzT} and \ref{app:The Error Term and Estimate} are taken into account. This result provides new insight, particularly with regard to the geometric interpretation of the terms of the formula; see Remark \ref{rem:geom-interp-of-zeta(0)}. Also, it introduces the proper conceptual framework and confirms that fractal sprays are clearly the higher-dimensional counterpart of fractal strings. In a similar vein, we will see from Theorem~\ref{thm:self-similar case} (the tube formula for self-similar tilings) that the self-similar tilings are the natural higher-dimensional analogue of self-similar fractal strings.

Although our primary goal in this paper is to obtain a tube formula for self-similar tilings, we state our main result for the more general class of fractal sprays, as we expect it to be useful in the study of other fractal structures and tilings to be investigated in future work. The key special case of self-similar tilings is stated in Theorem~\ref{thm:self-similar case} of \S\ref{ssec:The self-similar case}.

\subsection{Statement of the tube formula}
\label{ssec:statement of the tube formula}

We will prove the tube formula first for the more general case of fractal sprays, and then refine this result to obtain the formula for self-similar tilings. 

\begin{defn}
  \label{def:geometric zeta fn of a fractal spray}
  Let \tiling be a fractal spray with a single monophase generator \gen. Then the \emph{\tubezf zeta function} (or \emph{volume zeta function}) of \tiling is
  \begin{align}
    \label{eqn:monophase-geometric-zeta}
    \gzT(\ge,s) :=& \; \ge^{d-s} \gzs(s) \sum_{k=0}^d \frac{\genir^{s-k}}{s-k} \crv_{k}.
  \end{align}
  Here, $\crv_k = \crv_k(\gen)$ 
  as defined in Definition~\ref{def:monophase}.
\end{defn}

It turns out that \gzT is a meromorphic distribution-valued function for each fixed $s \in W$, where $W \ci \bC$ is the \emph{window} defined in Definition~\ref{def:screen}. This verification is given in Definition~\ref{def:d-valued mero} and Theorem~\ref{thm:gzT is a d-valued mero} of Appendix~\ref{app:Rigourous Definition of gzT}. Considered as a distribution, the action of $\gzT(s,\cdot)$ on a test function $\testfn \in C_c^\iy(0,\iy)$ is given by \begin{equation}
\label{eqn:def of action of gzT on gf}
  \left\la \gzT(\ge,s),\testfn(\ge)\right\ra
  = \int_0^\iy \gzT(\ge,s) \testfn(\ge) \,d\ge.
\end{equation}
Here, $C_c^\iy(0,\iy)$ is the space of smooth functions with compact support contained in $(0,\iy)$. At first glance, it may appear strange that something as concretely geometric as a tube formula is given distributionaly. However, the flexibility of the distributional framework allows the proof to proceed; see \cite[Rem.~5.20]{FGCD}.

\begin{remark}
  \label{rem:disc of the new gzT}
  The presentation here 
  differs slightly from that given in \cite{FGCD}, wherein the ``geometric zeta function'' is actually closer to what we call the scaling zeta function here. The general tube formula \eqref{eqn:1-dim tube formula} involves the one-dimensional case of \gzT, although this is not explicitly stated. For several reasons, it behooves one to think of \gzT as the zeta function most naturally associated with the geometric properties of the spray (or tiling), especially as pertains to the tube formula:
  \begin{Cond}
    \item The function \gzT arises naturally in the expression of the tube formula for the tiling, as will be seen in Theorem~\ref{thm:fractal-spray-tube-formula} and Theorem~\ref{thm:self-similar case}.
    \item It is the poles of $\gzT(\ge,s)$ that naturally index the sum appearing in $V_\tiling$, and the residues of \gzT that give the actual volume.
    \item Using \gzT leads to the natural unification of expressions which previously appeared unrelated; compare \eqref{eqn:main cor} to
      \eqref{eqn:single gen result} in Corollary~\ref{cor:self-similar-one-generator}.
  \end{Cond}
  Thus, the function \gzT encodes all the geometric information of \tiling as pertains to its tube formula.
  In Remark~\ref{rem:origin_of_mysterious_linear_term} we discuss how the unification mentioned in (ii) leads to a geometric interpretation of the term $\{2\ge \gzs(0)\}$ that appears in \eqref{eqn:1-dim tube formula}.
\end{remark}

\begin{defn}
  \label{def:visible complex dimensions of a fractal spray}
  The set of \emph{complex dimensions of a fractal spray} is
  \begin{align}
    \label{eqn:spray dimensions}
    \DT := \Ds 
              \cup \intDim,
  \end{align}
  where $\Ds$ is the set of poles of \gzs, as in \eqref{eqn:def:visible-scaling-dimns}.  When a window $W$ has been specified, the set of \emph{visible} complex dimensions is $\DT(W) := \DT \cap W$, and $\Ds(W) = \Ds \cap W$ is the set of visible scaling dimensions.
  Thus, $\DT(W)$ consists of the visible scaling dimensions and the visible ``integral dimensions'' of the spray. Furthermore, the poles of \gzT are all contained in $\DT$. Note that $\DT(W)$ is a discrete subset of $W \ci \bC$, and hence is countable. 
\end{defn}

\begin{theorem}[Tube formula for fractal sprays]
  \label{thm:fractal-spray-tube-formula}
  Let \tiling be a fractal spray on the monophase generator \gen, with generating inradius $\genir = \gr(\gen) > 0$, and scaling measure \ghs. Assume that \gzs is \slow on a screen $S$ which avoids the dimensions in $\DT(W)$.
  Then for test functions in $C_c^\iy(0,\iy)$,
  the $d$-dimensional volume of the inner tubular neighbourhood of the spray is
  given by the following distributional explicit formula:
  \begin{align}
    \label{eqn:fractal-spray-main-result}
    V(\tiling,\ge) &=   \sum_{\gw \in \DT(W)}
    \res{\gzT(\ge,s)} + \sR(\ge),
  \end{align}
  where the sum ranges over the set \eqref{eqn:spray dimensions}
  of visible integral and scaling dimensions of the spray.
  Here, the error term $\sR(\ge)$ is
  given by
  \begin{align}
    \label{def:main error term}
    \sR(\ge) 
    = \frac1{2\gp \ii}\int_S \gzT(\ge,s) \, ds,
  \end{align}
  and estimated by
  \begin{align}
    \label{eqn:main error est}
    \sR(\ge) &=  O(\ge^{d-\sup S}), \qq \text{as } \ge \to 0^+.
  \end{align}
\end{theorem}
In the case that $\gw \in \Ds(W) \cap \intDim$, then the corresponding term $\res{\gzT(\ge,s)}$ appears only once in the sum in \eqref{eqn:fractal-spray-main-result}. 
As a distributional formula, \eqref{eqn:fractal-spray-main-result} is valid when applied to test functions $\testfn \in C_c^\iy(0,\iy)$. The order of the distributional error term as in \eqref{eqn:main error est} is defined in Definition~\ref{def:order of a dist}. There is a version of this theorem in which the error term vanishes identically; it is presented in Corollary~\ref{cor:strongly languid fractal spray tube formula}. Also, the special case of self-similar tilings is presented in Theorem~\ref{thm:self-similar case}. The following proof relies heavily on the material in \S\ref{sec:The distributional theory of fractal strings}; the reader may wish to review this material before proceeding.

\subsection{Proof of the tube formula}
\label{ssec:proof of the tube formula}

The reader may now wish to review \S\ref{sec:The distributional theory of fractal strings} before proceeding, as the proof uses these explicit formulas and distributional techniques from \cite{FGCD}.

\begin{proof}[Proof of Theorem~\ref{thm:fractal-spray-tube-formula}]
Recall that we view $V(\tiling,\ge)$ as a
distribution,\footnote{Indeed, $V(\tiling,\ge)$ is clearly continuous and
bounded (by the total volume of the spray); hence it defines a
locally integrable function on $(0,\iy)$.} so we understand
$V(\tiling,\ge) = \la\ghs, \gtf \ra$ by computing its action on
a test function \testfn:
\begin{align}
  \left \la V(\tiling,\ge), \testfn \right \ra =
  \left \la \la \ghs, \gtf \ra, \testfn \right \ra
  &= \int_0^\iy \left(\int_0^\iy \gtf(x,\ge) d\ghs(x)\right) \testfn(\ge) \, d\ge \notag \\
  &= \int_0^\iy \int_0^\iy \gtf(x,\ge) \testfn(\ge) \, d\ge \, d\ghs(x) \notag \\
  &= \left\la \ghs, \la \gtf, \testfn \ra \right\ra.
  \label{eqn:action of V as action of eta}
\end{align}
Now, writing $\crv_k = \crv_k(\gen)$, we use
\eqref{eqn:thm:monophase-formula} to compute
\begin{align}
  \la \gtf, \testfn \ra
  &= \int_0^\iy \gtf(x,\ge) \testfn(\ge) \, d\ge \notag \\
  &= \int_0^\iy \sum_{k=0}^{d} \crv_{k}(\tfrac1x \gen,\ge) \ge^{d-k} \testfn(\ge) \,d\ge \notag \\
  &= \sum_{k=0}^{d-1} \int_0^\iy \crv_{k} \gc_{[0,\genir/x)}(\ge) x^{-k} \ge^{d-k} \testfn(\ge) \,d\ge
    - \int_0^\iy \crv_{d} \gc_{[\genir/x,\iy)}(\ge) x^{-d} \testfn(\ge) \, d\ge \notag \\
  &= \sum_{k=0}^{d-1} \crv_{k} x^{-k} \int_0^{\genir/x} \ge^{d-k} \testfn(\ge) \,d\ge
    - \crv_{d} x^{-d} \int_{\genir/x}^\iy \testfn(\ge) \, d\ge \notag \\
  &= \sum_{k=0}^d \gf_{k}(x),
  \label{eqn:action of ve on phi}
\end{align}
where, for $x>0$, we have introduced
\begin{equation}
  \label{eqn:def-of-gf_i}
  \gf_{k}(x) :=
  \begin{cases}
    \crv_{k} x^{-k} \int_0^{\genir/x} \ge^{d-k} \testfn(\ge) \,d\ge, & 0 \leq k \leq d-1, \\
    \crv_{k} x^{-k} \int_\iy^{\genir/x} \testfn(\ge) \, d\ge, &k=d,
  \end{cases}
\end{equation}
in the last line. Caution: $\gf_{k}$ is a function of $x$, whereas
$\testfn$ is a function of \ge. 

  Putting (\ref{eqn:action of ve on
phi}) into (\ref{eqn:action of V as action of eta}), we obtain
\begin{align}
  \left \la V(\tiling,\cdot), \testfn \right \ra
  &= \left\la \ghs, \sum_{k=0}^d \gf_{k} \right\ra
   = \sum_{k=0}^d \left\la \ghs, \gf_{k} \right\ra.
  \label{eqn:action of V first sum}
\end{align}
To apply Theorem~\ref{thm:ext dist explicit formula}, we must
first check that the functions $\gf_{k}$ satisfy the
hypotheses \eqref{eqn:test fn first req for extended
formula}--\eqref{eqn:test fn second req for extended formula}.
Recall that $\testfn \in C_c^\iy(0,\iy)$.

For $k<d$, \eqref{eqn:test fn first req for extended formula}
is satisfied because for large $x$, the corresponding integral in
\eqref{eqn:def-of-gf_i} is taken over a set outside the
(compact) support of \testfn. This gives $\gf_{k}(x) = 0$ for
sufficiently large $x$, and it is clear that, \emph{a
fortiori}, the \nth derivative of $\gf_k$ satisfies
\begin{equation}
  \label{eqn:gf_i satisfies first cond}
  \gf_{k}^{(n)}(x) = O(x^{-n-D-\gd}) \q \text{for } x \to \iy. 
\end{equation}

To see that \eqref{eqn:test fn second req for extended formula}
is satisfied, note that $\testfn$ vanishes for $x$ sufficiently
large and thus 
\begin{align*}
  \gf_{k}(x)
  &= \crv_{k} x^{-k} \int_0^{\iy} \ge^{d-k} \testfn(\ge) \,d\ge
  \q \text{for } x \approx 0,
\end{align*}
i.e., $\gf_{k}(x) = a_{k} x^{-k}$ for all small enough
$x>0$, where $a_{k}$ is the constant
\begin{equation}
  \label{def:ai}
  a_{k} := \crv_{k} \int_0^{\iy} \ge^{d-k} \testfn(\ge) \,d\ge
  = \crv_{k} \testfnt(d-i+1)
  = \lim_{x \to 0^+} x^k\gf_{k}(x).
\end{equation}
Here \testfnt is the Mellin transform of \testfn, as in
\eqref{def:Mellin transform}.

Thus, the expansion \eqref{eqn:test fn second req for extended formula} for the test function $\gf_{k}$ consists of only one term, and for each $n=0,1,\dots,N$,\footnote{Recall that \ghs is languid of order \languidOrder and that $N = [\languidOrder]+2$ in the hypotheses of Theorem~\ref{thm:ext dist explicit formula}. }
\begin{equation}
  \label{eqn:gf_i satisfies second cond}
  \gf_{k}^{(n)}(x) = \tfrac{d^n}{dx^n}\left[a_{k} x^{-k}\right] = O(x^{-n-k})
  \q \text{for } x \to 0^+. 
\end{equation}
A key point is that since \testfn is smooth, \eqref{eqn:gf_i satisfies first cond} and \eqref{eqn:gf_i satisfies second cond} will hold for each $n=0,1,\dots,N$, as required by Theorem~\ref{thm:ext dist explicit formula}. Since the expansion of $\gf_{k}$ has only one term, the only \ga in the sum is $\ga=k$. Thus $a_{k}$ is the constant corresponding to $a_\ga$ in \eqref{eqn:test fn second req for extended formula}.

Applying Theorem~\ref{thm:ext dist explicit formula} in the case when $k<d$, \eqref{eqn:action of dist on test fn} becomes
\begin{align}
  \label{eqn:gh action for phij}
    \left\la \ghs, \gf_{k}\right\ra
    &= \sum_{\gw \in \Ds(W)} \negsp[10] \res{\gzs(s)\gfti(s)}
       \; + \; \{a_{k} \gzs(k)\}_{k \in W \less \Ds} \notag \\
       &\hstr[12] \; + \; \tfrac1{2\gp \ii} \int_S \gzs(s) \gfti(s)\,ds,
\end{align}
where the term in braces is to be included iff $k \in W \less \Ds$. Here and henceforth, \gfti denotes the Mellin transform of \gfi given by
\begin{equation}
  \label{def:BjqMellin}
  \gfti(s) = \int_0^{\iy} x^{s-1} \gf_k(x) \,dx.
\end{equation}

The case when $k=d$ is similar (or antisimilar). The compact
support of \testfn again gives
\begin{align}
  \gf_{d}(x)
  &= \crv_{d} x^{-d} \int_\iy^0 \testfn(\ge) \, d\ge,
  \q\text{for } x \to \iy,
\end{align}
so that for some positive constant $c$, and for all
sufficiently large $x$, we have $\crv_{d}(x) = cx^{-d}$. Hence
\begin{equation}
  \label{eqn:gf_d satisfies first cond}
  \gf_{d}^{(n)}(x)
  = O(x^{-n-d}) \q \text{for } x \to \iy, \, \forall n \geq 0,
\end{equation}
and \eqref{eqn:test fn first req for extended formula} is
satisfied. For very small $x$, the integral in the definition of
$\gf_{d}(x)$ is taken over an interval outside the support of
\testfn, and hence $\gf_{d}(x)=0$ for $x \approx 0$. Then
clearly \eqref{eqn:test fn second req for extended formula} is
satisfied:
\begin{equation}
  \label{eqn:gf_d satisfies second cond}
  \gf_{d}^{(n)}(x) = 0 \q \text{for } x \to 0^+, \, \forall n \geq 0.
\end{equation}
An immediate consequence of \eqref{eqn:gf_d satisfies second
cond} is that for $k=d$ in \eqref{def:ai}, the constant term
is
\begin{equation}
  \label{eqn:ad=0}
  a_d = \lim_{x \to 0} x^d \gf_d(x) = 0,
\end{equation}
and compared with \eqref{eqn:gh action for phij} we have one term
less in
\begin{equation}
  \label{eqn:gh action for phid}
  \left\la \ghs, \gf_d \right\ra
    = \sum_{\gw \in \Ds(W)} \negsp[10] \res{\gzs(s)\gft_{d}(s)}
    + \tfrac1{2\gp \ii} \int_S \gzs(s) \gft_{d}(s)\,ds.
\end{equation}

As in \eqref{def:BjqMellin}, denote the Mellin transform of the
function \gy by $\widetilde \gy$ and compute
\begin{align}
  \gfti(s)
   = \int_0^\iy x^{s-1} \gf_{k}(x) \, dx
  &= \crv_{k} \int_0^\iy x^{s-k-1} \int_0^{\genir/x} \ge^{d-k} \testfn(\ge) \,d\ge \, dx \notag \\
  &= \crv_{k} \int_0^\iy \left(\int_0^{\genir/\ge} x^{s-k-1} \, dx \right) \ge^{d-k} \testfn(\ge) \,d\ge \notag \\
  &= \frac{\crv_{k}}{s-k} \int_0^\iy \genir^{s-k} \ge^{k-s} \ge^{d-k} \testfn(\ge) \, d\ge \notag \\
  &= \genir^{s-k} \frac{\crv_{k}}{s-k} \testfnt(d-s+1).
  \label{eqn:gfti of s}
\end{align}
By a similar calculation,
\begin{align}
  \gft_{d}(s)
   &= \genir^{s-d} \frac{\crv_{d}}{s-d} \testfnt(d-s+1).
  \label{eqn:gftd of s}
\end{align}
%
Note that for $0 \leq k < d-1$, \eqref{eqn:gfti of s} is valid
for $\Re s > k$, and for $k=d$, \eqref{eqn:gftd of s} is
valid for $\Re s < k$. Thus both are valid in the strip $d-1 <
\Re s < d$, and hence by analytic (meromorphic) continuation,
they are valid everywhere in \bC. Indeed, by
Corollary~\ref{cor:holomorphicity of gft}, \testfnt is entire.

We return to the evaluation of \eqref{eqn:action of V first
sum}, applying Theorem~\ref{thm:ext dist explicit formula} to find
the action of \ghs on the test function $\gf_{k}$, for
$k=0\dots,d$. Substituting (\ref{eqn:gfti of s}) and
\eqref{eqn:gftd of s} into (\ref{eqn:gh action for phij}) gives
\begin{align}
  \label{eqn:ghj action as Mellins}
  \left\la \ghs, \gf_{k}\right\ra
    &= \sum_{\gw \in \Ds(W)} \negsp[10] \res{\gzs(s) \frac{\genir^{s-k} \crv_k}{s-k} \testfnt(d-s+1)} \\
    &\hstr[8] \; + \; \{a_{k} \gzs(k)\}_{k \in W \less \Ds}
    \; + \; \la\sR_{k},\testfn\ra, \notag
\end{align}
where $\sR_{k}$ is defined by
\begin{align}
  \label{def:partial error Rq}
  \la\sR_{k}, \testfn\ra := \tfrac1{2\gp \ii} \int_S \gzs(s) \gfti(s)\,ds.
\end{align}
Substituting \eqref{eqn:ghj action as Mellins} into
\eqref{eqn:action of V first sum}, we obtain
\begin{align}
  \left \la V(\tiling,\ge), \testfn \right \ra
  &= \sum_{k=0}^d
    \sum_{\gw \in \Ds(W)} \negsp[10] \res{\gzs(s) \frac{\genir^{s-k} \crv_k}{s-k} \testfnt(d-s+1)}
     \notag \\ &
    \hstr[4] + \sum_{k=0}^d \{a_{k} \gzs(k)\}_{k \in W \less \Ds}
    + \sum_{k=0}^d \la \sR_{k}(\ge),\testfn(\ge)\ra.
    \label{eqn:V(e)-as-dist}
\end{align}
Recall from \eqref{eqn:ad=0} that the \nth[d] term is \(a_{d} =
0\), so the top term of the second sum vanishes. Note that at
each such $k$ we have a residue
\begin{align}
  \label{eqn:res at j}
  \res[k]{\gzs(s) \frac{\genir^{s-k} \crv_{k}}{s-k} \testfnt(d-s+1)}
  &= \crv_{k} \lim_{s \to k} \gzs(s) \genir^{s-k} \testfnt(d-s+1) \notag \\
  &= \crv_{k} \gzs(k) \testfnt(d-i+1) \notag \vstr \\
  &= a_k \gzs(k).\vstr
\end{align}
Since the terms of the second sum of \eqref{eqn:V(e)-as-dist}
are only included for $k \in W \less \Ds(W)$,
we can use \eqref{eqn:spray dimensions} and \eqref{eqn:res at
j} to combine the last two sums of \eqref{eqn:V(e)-as-dist}
without losing or duplicating terms:
\begin{align*}
  \left \la V(\tiling,\ge), \testfn \right \ra
  &= \sum_{\gw \in \DT(W)} \negsp[10] \res{\testfnt(d-s+1) \gzs(s)
     \sum_{k=0}^d \frac{\genir^{s-k} \crv_k}{s-k} }
    + \left\la \sR(\ge),\testfn(\ge) \right\ra,
\end{align*}
where $\sR(\ge) := \sum_{k=0}^d \sR_{k}(\ge)$. This may also
be written as the distribution
\begin{align}
  \label{eqn:V as a function}
   V(\tiling,\ge) &=  \sum_{\gw \in \Ds(W)} \negsp[10]
    \res{\ge^{d-s} \gzs(s) \sum_{k=0}^d \frac{\genir^{s-k} \crv_k }{s-k}}
     + \sR(\ge).
\end{align}

  This completes the proof of \eqref{eqn:fractal-spray-main-result}. All that
  remains is the verification of the expression \eqref{def:main error
  term} for the error term, and error estimate \eqref{eqn:main error
  est}. Due to their technical and specialized nature, we leave the
  proofs of \eqref{def:main error term} and \eqref{eqn:main error est}
  to Appendix~\ref{app:The Error Term and Estimate}.
\end{proof}

\section{Extensions and Consequences: the Tube Formula for Self-Similar Tilings}
\label{sec:Extensions and Consequences}


Recall from \S\ref{ssec:Tilings with multiple generators} that the results of \S\ref{sec:Extensions and Consequences} may easily be extended to multiple generators simply by taking the corresponding finite sum.
The next corollary indicates that when \gzs is strongly languid, one may take $W=\bC$ in the previous theorem and the error term will vanish identically.

\begin{cor}[Tube formula for strongly languid fractal sprays]
  \label{cor:strongly languid fractal spray tube formula}
  Let \tiling be a fractal spray on the monophase generator \gen with \slow scaling measure \ghs, and additionally assume that \gzs is \emph{strongly} languid, and hence that $W=\bC$.
  Then
  \begin{align}
    \label{eqn:fractal-spray-L2-result}
    V(\tiling,\ge) &=  \sum_{\gw \in \DT} \res{\gzT(\ge,s)},
  \end{align}
  where $\DT = \DT(\bC)$ is the set of complex dimensions of \tiling, as in \eqref{eqn:spray dimensions}.
\end{cor}
\begin{proof}
  This is immediate upon combining Theorem~\ref{thm:ext_distributional_without error}
  (the extended distributional formula without error term) with
  the proof of Theorem~\ref{thm:fractal-spray-tube-formula}.
  One finds that $\R_{k} \equiv 0$ for each $k=0,1,\dots,d$ in \eqref{def:partial error Rq} and thus $\R \equiv 0$ in \eqref{eqn:V as a function}.
\end{proof}

\begin{remark}[Reality principle]
  \label{rem:reality principle}
  The nonreal complex dimensions appear in complex conjugate pairs and
  produce terms with coefficients which are also complex
  conjugates, in the general tube formula for fractal sprays. This
  ensures that formulas \eqref{eqn:fractal-spray-main-result} and \eqref{eqn:fractal-spray-L2-result}--\eqref{eqn:main cor} are real-valued.
\end{remark}

\subsection{The self-similar case}
  \label{ssec:The self-similar case}

Self-similar strings automatically satisfy the more stringent hypothesis of being
\emph{strongly languid}, as in Definition~\ref{def:strongly languid}. This automatically entails that Corollary~\ref{cor:strongly languid fractal spray tube formula} holds,\footnote{This is essentially because Theorem~\ref{thm:ext dist explicit formula} and Theorem~\ref{thm:tube formula for fractal strings} hold without error term. This is discussed further in \cite[Thm.~5.27]{FGCD}, and the end of \cite[Thm.~8.1]{FGCD}. A general discussion of the strongly languid case may be found in \cite[Def.~5.3]{FGCD}, and an argument showing that all self-similar strings are strongly languid is given in \cite[\S6.4]{FGCD}.} so the window may be taken to be all of \bC and the error term vanishes identically, i.e., $\sR(\ge) \equiv 0$. Hence Theorem~\ref{thm:fractal-spray-tube-formula} may be strengthened for self-similar tilings as in Theorem~\ref{thm:self-similar case}.

\begin{theorem}[Tube formula for self-similar tilings]
  \label{thm:self-similar case}
  Let
  \(\tiling = \{\simt_w \gen\},\)
  be a self-similar tiling with monophase generator \gen and
  tubular zeta function $\gzT$.
  Then the \mbox{$d$-dimensional} volume of the inner tubular neighbourhood of \tiling
  is given by the following distributional explicit formula:
  \begin{align}
    \label{eqn:main cor}
    V(\tiling,\ge) &=  \sum_{\gw \in \DT} \negsp[1] \res{\gzT(\ge,s)},
  \end{align}
  where $\DT = \DT(\bC) = \Ds(\bC) \cup \intDim$ is the set of complex dimensions of \tiling.
\end{theorem}
  \begin{proof}
    The proof follows \cite[\S6.4]{FGCD}.
    According to Theorem~\ref{thm:scaling-zeta-fn-simplified},
    the scaling zeta function of a self-similar tiling has the form
    \[\gzs(s) = \frac{1}{1 - \sum_{\j=1}^J r_\j^s}.\]
    Let $r_J$ be the smallest scaling ratio. Then from
    \[|\gzs(s)| \ll \left(\frac{1}{r_J}\right)^{-|\gs|}
      \q \text{as } \gs=\Re(s) \to -\iy,\]
    we deduce that \gzT is strongly languid and therefore apply
    Corollary~\ref{cor:strongly languid fractal spray tube formula}.
    This argument follows from the analogous ideas regarding
    self-similar strings, which may be found in \cite[\S8.4]{FGCD}.
  \end{proof}

\begin{remark}
  Theorem~\ref{thm:self-similar case} provides a higher-dimensional counterpart of the tube formula obtained for self-similar strings in \cite[\S8.4]{FGCD}. It should be noted that Theorem~\ref{thm:self-similar case} applies to a slightly smaller class of test functions than Theorem~\ref{thm:fractal-spray-tube-formula}. Indeed, the support of the test functions must be bounded away from 0 by $\gm_d(\hull) \genir/r_J$, where $\hull=[\attr]$ is the hull of the attractor (as in \S\ref{sec:The self-similar tiling}), $\genir$ is the smallest generating inradius (as in \eqref{eqn:genir ordering}), and $r_J$ is the smallest scaling ratio of \simt (as in \eqref{eqn:scaling ratio ordering}). This technicality is discussed further in \cite{FGCD}, Def.~5.3 and Thm.~5.27, \S6.4, and Thm.~8.1.
\end{remark}


\begin{cor}[Measurability and the lattice/nonlattice dichotomy]
  \label{cor:meas and dichotomy}
  Under mild conditions on the residues of \gzT, a self-similar tiling is Minkowski measurable if and only if it is nonlattice.
\end{cor}
  \begin{proof}[Sketch of proof]
    We define a self-similar tiling \tiling to be Minkowski measurable iff
    \begin{align}
      \label{def:Mink meas}
      0 < \lim_{\ge \to 0^+} V(\tiling,\ge) \ge^{-(d-D)} < \iy,
    \end{align}
    i.e., if the limit in \eqref{def:Mink meas} exists \emph{and} takes a value in $(0,\iy)$. A tiling has infinitely many complex dimensions with real part $D$ iff it is lattice type, as mentioned in \S\ref{ssec:comparison to FGNT}. Furthermore, all the poles with real part $D$ are simple in that case. A glance at \eqref{eqn:single-gen-conceptual-result} then shows that $V(\tiling,\ge) \ge^{-(d-D)}$ is a sum containing infinitely many purely oscillatory terms $c_\gw \ge^{\ii n\per}$, $n \in \bZ$, where \per is some fixed period. Thus, the limit \eqref{def:Mink meas} cannot exist; see also \cite[\S8.4.2]{FGCD}.\footnote{It is shown in \cite{FGCD} that infinitely many coefficients $c_\gw$ are nonzero for $\Re \gw = D$.}  Conversely, the tiling is nonlattice iff $D$ is the only complex dimension with real part $D$. In this case, $D$ is simple and no term in the sum $V(\tiling,\ge) \ge^{-(d-D)}$ is purely oscillatory; thus the tiling \tiling is measurable. See also \cite[\S8.4.4]{FGCD}.
  \end{proof}

\emph{Note added in proof:} Please see Remark~\ref{rem:note-added-in-proof} for the ``mild conditions'' mentioned in the statement of Corollary~\ref{cor:meas and dichotomy}, and the scope of Remark~\ref{rem:dimensions-and-measurability} and Remark~\ref{rem:combine-FGCD-with-this}.

\begin{remark}
  \label{rem:dimensions-and-measurability}
  In \cite[\S8.3--8.4]{FGCD}, it is shown that a self-similar fractal string (i.e., a \linebreak \mbox{1-dimensional} self-similar tiling) is Minkowski measurable if and only if it is nonlattice. Gatzouras showed in \cite{Gat} that nonlattice self-similar subsets of \bRd are Minkowski measurable, thereby extending to higher dimensions a result in \cite{La3}, \cite{Fal2} and partially proving the geometric part of \cite[Conj.~3]{La3}. The previous result gives a complete characterization of self-similar tilings in \bRd as nonlattice if and only if they are Minkowski measurable. With the exception of Remark~\ref{rem:nonlattice Koch}, each of the examples discussed in \S\ref{sec:Tube formula examples} is lattice and hence not Minkowski measurable. Our results, however, apply to nonlattice tilings as well. A more detailed proof of Corollary~\ref{cor:meas and dichotomy} is possible via truncation, by using the screen and window technique of \cite[Thm.~5.31 and Thm.~8.36]{FGCD}.
\end{remark}

The following corollary of Theorem~\ref{thm:self-similar case}
will be used in \S\ref{sec:Tube formula examples}.
\begin{cor}
  \label{cor:self-similar-one-generator}
  If, in addition to the hypotheses of Theorem~\ref{thm:self-similar case},
  $\gzT(s)$ has only simple poles, then
  \begin{align}
    \label{eqn:single gen result}
    V(\tiling,\ge)
    &= \sum_{\gw \in \Ds} \res{\gzs(s)} \ge^{d-\gw} \sum_{k=0}^{d} 
      \tfrac{\genir^{\gw-k}}{\gw-k} \crv_k
      + \sum_{k=0}^{d-1} \crv_k \gzs(k) \ge^{d-k}.
  \end{align}
\end{cor}

It is not an error that the first sum extends to $d$ in
\eqref{eqn:single gen result}, while the second stops at $d-1$;
see \eqref{eqn:ad=0}. Note that in Corollary~\ref{cor:self-similar-one-generator}, \Ds does not contain any integer
$k=0,1,\dots,d-1$, because this would imply that \gzT has a pole
of multiplicity at least 2 at such an integer. In general, at most one integer can
possibly be a pole of \gzs; see \S\ref{ssec:comparison to
FGNT}.

\begin{remark}
  \label{rem:conceptual version of the tube formula}
  For self-similar tilings satisfying the hypotheses of
  Corollary~\ref{cor:self-similar-one-generator}, it is clear
  that the general form of the tube formula is
  \begin{align}
    \label{eqn:single-gen-conceptual-result}
    V(\tiling,\ge)
    &= \sum_{\gw \in \DT}  c_\gw \ge^{d-\gw},
  \end{align}
  where for each fixed $\gw \in \Ds$,
  \begin{align}
    \label{def:c_wj}
    c_\gw
    &:=  \res{\gzs(s)} \sum_{k=0}^d \frac{\genir^{\gw-k}}{\gw-k} \crv_k.
  \end{align}
  Note that when $\gw = k \in \{0,1,\dots,d-1\}$, one has
  $c_\gw = c_k = \gzs(k) \crv_k$.

%
\end{remark}

\begin{remark}
  The oscillatory nature of the geometry of \tiling is apparent in \eqref{eqn:single-gen-conceptual-result}. In particular, the existence of the limit in \eqref{def:Mink meas} can be determined in the nonlattice case by examining \eqref{eqn:single-gen-conceptual-result} and \DT.
\end{remark}

\begin{remark}
  \label{rem:gaps vs generators}
  In the literature regarding the 1-dimensional case \cite{La-vF3}, \cite{FGCD}, \cite{Fra}, the terms ``gaps'' and ``multiple gaps'' have been used where we have used ``generators''.
\end{remark}

\begin{remark}[Comparison of $V(\tiling,\ge)$ with the Steiner formula]
  \label{rem:comparison with Steiner}
  In the trivial situation when the spray consists only of finitely
  many scaled copies of a monophase generator (so the scaling measure
  \ghs is supported on a finite set), the zeta function \gzs will have
  no poles in \bC.
  Therefore, the tube formula becomes a sum
  over only the numbers $0,1,\dots,d-1$ (recall from \eqref{eqn:ad=0}
  that $a_{d}=0$, so the \nth[d] summand vanishes), for which the residues
  simplify greatly as in \eqref{eqn:res at j}.
%
  In this case, $\gzs(k) = \gr_1^k + \dots + \gr_J^k$, so each residue
  from \eqref{eqn:res at j}
  becomes a finite sum
  \begin{align*}
    \gzs(k)\crv_{k}(\ge)
    &= \gr_1^k \crv_{k} \ge^{d-k} + \dots + \gr_J^k \crv_{k} \ge^{d-k} \\
    &= \crv_{k} (r_{w_1} \gen) \ge^{d-k}
    + \dots + \crv_{k} (r_{w_J} \gen) \ge^{d-k},
  \end{align*}
  where $J$ is the number of scaled copies of the generator
  \gen, and $r_{w_j}$ is the corresponding scaling factor. Thus,
  for each $j=1,\dots,J$, we obtain a monophase formula for the scaled
  basic shape $r_{w_j} \gen$. 
\end{remark}


\begin{remark}\label{rem:combine-FGCD-with-this}[Combining the results of this paper with \cite{FGCD}.]
  Since, as was noted in \S\ref{ssec:comparison to FGNT}, the structure of the scaling complex dimensions of a self-similar tiling in \bRd is the same as in the $1$-dimensional case, we could state an analogue of each of the theorems given in \cite[\S8.4]{FGCD}, whether in the lattice case (\cite[\S8.4.2]{FGCD}) or in the nonlattice case (\cite[\S8.4.4]{FGCD}). In particular, we can apply Theorem~\ref{thm:fractal-spray-tube-formula} with a suitable window $W$ (and use the Diophantine approximation techniques of \cite[Ch.~3]{FGCD}) in order to obtain the exact higher-dimensional analogues of the tube formulas with error term stated in \cite[Cor.~8.27]{FGCD} and \cite[Thm.~8.37 and Eqn.~(8.71)]{FGCD}, in the lattice and nonlattice case, respectively. Finally,in the lattice case, following \cite[\S8.4.3]{FGCD}, even though the self-similar tiling is not Minkowski measurable, we could calculate its `average Minkowski content' (cf.~\cite[Def.~8.29]{FGCD} with 1 replaced by $d$) and obtain the $d$-dimensional analogue of \cite[Thm.~8.30]{FGCD}. In order to avoid redundancies, we will avoid formulating explicitly any of these consequences of our results in this paper. An attentive reader of \cite{FGCD} should easily be able to combine that material with the present results and obtain such useful corollaries. 
\end{remark}

\subsection{Recovering the tube formula for
fractal strings} \label{ssec:The tube formula for fractal
strings}

In this section, we discuss the \mbox{1-dimensional} tube formula of Theorem~\ref{thm:tube formula for fractal strings} which is true for general (i.e., not necessarily self-similar) fractal strings and which can be recovered from Theorem~\ref{thm:fractal-spray-tube-formula}. Suppose $\sL=\{\ell_n\}_{n=1}^\iy$ is a languid fractal string with associated measure $\gh_\sL=\sum_{n=1}^\iy \gd_{1/\ell_n}$, as in \eqref{def:ghL}, and associated zeta function $\gzL = \sum_{n=1}^\iy \ell_n^s$, as in \eqref{def:string zeta fn}. Considering the string now as a tiling, write \sL as $L = \{L_n\}_{n=1}^\iy$ to emphasize the fact that we are thinking of it as a spray instead of as a string. Take the spray $L$ to have as its single generator the interval $\gen=(0,2)$, so that $L$ has inradii $\gr_n = \tfrac12\ell_n$ and that the length of $\gr_n \gen$ is $\ell_n$. Now the scaling measure is $\ghs 
= \sum_{n=1}^\iy \gd_{2/\ell_n}$ and the scaling zeta function is
\begin{align} \label{eqn:string-scaling-zeta-function}
  \gzs(s)
  = \sum_{n=1}^\iy \left(\frac{\ell_n}{2}\right)^s
  = 2^{-s} \gzL(s).
\end{align}
The generator is clearly monophase with $\crv_0 = 2$ and $\crv_1
= -2\genir$:
\begin{align}
  \gtf(x,\ge)
    = \begin{cases}
        2\ge, &\ge \leq \genir/x, \\
        2\genir/x, &\ge \geq \genir/x.
      \end{cases}
    \label{eqn:general string gtf}
\end{align}
One obtains the \tubezf zeta function of the (1-dimensional) tiling $L$ as
\begin{align}
  \label{eqn:general string inner product}
  \gz_L(\ge,s)
  &= \ge^{1-s} \gzs(s) \sum_{k=0}^1 \frac{\crv_k}{s-k}
   = \ge^{1-s} 2^{-s} \gzs(s) \left(\frac{2}{s} - \frac{2}{s-1}\right)
   =\frac{\gzL(s) (2\ge)^{1-s}}{s(1-s)},
\end{align}
by substituting in \eqref{eqn:string-scaling-zeta-function} in the last step.  
Then
from Theorem~\ref{thm:fractal-spray-tube-formula} we exactly recover the tube formula $V(L,\ge) = V(\sL,\ge)$ (and its error term) as given by Theorem~\ref{thm:tube formula for fractal strings}. Note that $\Ds = \sD_\sL$ by \eqref{eqn:string-scaling-zeta-function}. 


\begin{remark}
  \label{rem:origin_of_mysterious_linear_term}
  In addition to recovering a previously known formula, we also gain
  a geometric interpretation of the terms appearing in the 1-dimensional
  tube formula \eqref{eqn:1-dim tube formula}, in view
  of the previous computation.
  In particular, one sees that the linear term $\{2\ge \gzs(0)\}$
  has a geometric interpretation in terms of the inner Steiner formula
  for an interval, and can be dissected as
  \begin{equation}
    \label{eqn:dissection of 2ge gz(0)}
    2\ge \gzs(0) = \crv_0(\gen) \ge^{1-0} \gzs(0) = (-2) \gm_k(\gen) \ge^{d-k} \gzs(k),
  \end{equation}
  where $k=0$ and $d=1$.
  Note that $\gm_0(\gen)=-1$ is the Euler characteristic of an open
  interval.
  %
  This should be discussed further in \cite{FCM}.
\end{remark}


\section{Tube Formula Examples}
\label{sec:Tube formula examples}

Although Remark~\ref{rem:nonlattice Koch} discusses how one may
construct nonlattice examples, the other examples chosen in this
section are \emph{lattice} self-similar tilings, in the sense of
\S\ref{ssec:comparison to FGNT}. 
Also,
all examples in this section have monophase generators in the
sense of Definition~\ref{def:monophase}, as is verified in each case. The pentagasket of Example~\ref{exm:Pentagasket Tiling} is the only example given here of a self-similar tiling with multiple generators.

Moreover, the scaling zeta function \gzs of each example has only simple
poles, with a single line of complex dimensions distributed
periodically on the line $\Re s = D$. Thus, the tube formula may be
substantially simplified via Corollary~\ref{cor:self-similar-one-generator}.

\subsection{The Cantor tiling} 
\label{exm:Cantor Tiling}

The Cantor tiling \sC (called the Cantor string in \cite{FGCD}, \S1.1.2 and
\S2.3.1) is constructed via the self-similar system
\[\simt_1(x) = \tfrac x3, \qq
  \simt_2(x) = \tfrac{x+2}3.\]
The associated self-similar set \attr is the classical ternary
Cantor set, so $d=1$ and we have one scaling ratio \(r=\tfrac13,\)
and one generator $\gen=\left(\tfrac13,\tfrac23\right)$ with
generating inradius \(\genir=\tfrac16.\) The corresponding
self-similar string has inradii $\gr_m = \genir r^m$ with
multiplicity $2^m$, $m=0,1,2,\dots$, so the scaling zeta function is
\begin{align}
  \label{eqn:Cantor zetafn}
  \gzs(s) = \frac1{1 - 2 \cdot 3^{-s}},
\end{align}
and the scaling complex dimensions are
\begin{align}
  \label{eqn:Cantor dimensions}
  \Ds = \{D + \ii n\per \suth n \in \bZ\}
  \qq \text{for } D=\log_32, \; \per=\tfrac{2\gp}{\log3}.
\end{align}

We note that $\gzs(0)=-1$ and apply \eqref{eqn:general string inner
product} from the previous section to recover the following tube
formula for \sC (as obtained in \cite{FGCD}, \S1.1.2):
\begin{align}
  V(\sC,\ge)
  &= \frac1{2\log3} \sum_{n \in \bZ}
       \frac{(2\ge)^{1-D-\ii n\per}}{(D+\ii n\per)(1-D-\ii n\per)}
          - 2\ge.
  \label{eqn:Cantor final by formula}
\end{align}

Alternatively, this may be written as a series in
$\left(\frac{\ge}{\genir}\right)$ as
\begin{align}
  V(\sC,\ge)
  &= \frac1{3\log3} \sum_{n \in \bZ}
     \left(\frac{1}{D+\ii n\per}-\frac1{D-1+\ii n\per}\right)\left(\frac{\ge}{\genir}\right)^{1-D-\ii n\per} -
     2\ge,
  \label{eqn:Cantor final in 6eps}
\end{align}
with $\genir=\tfrac16$, $D=\log_32$, and $\per=2\gp/\log3$. It is
this form of the tube formula which is closer in appearance to the
following examples.

\subsection{The Koch tiling}
\label{ssec:The Koch Tiling}

\begin{figure}
  \scalebox{0.80}{\includegraphics{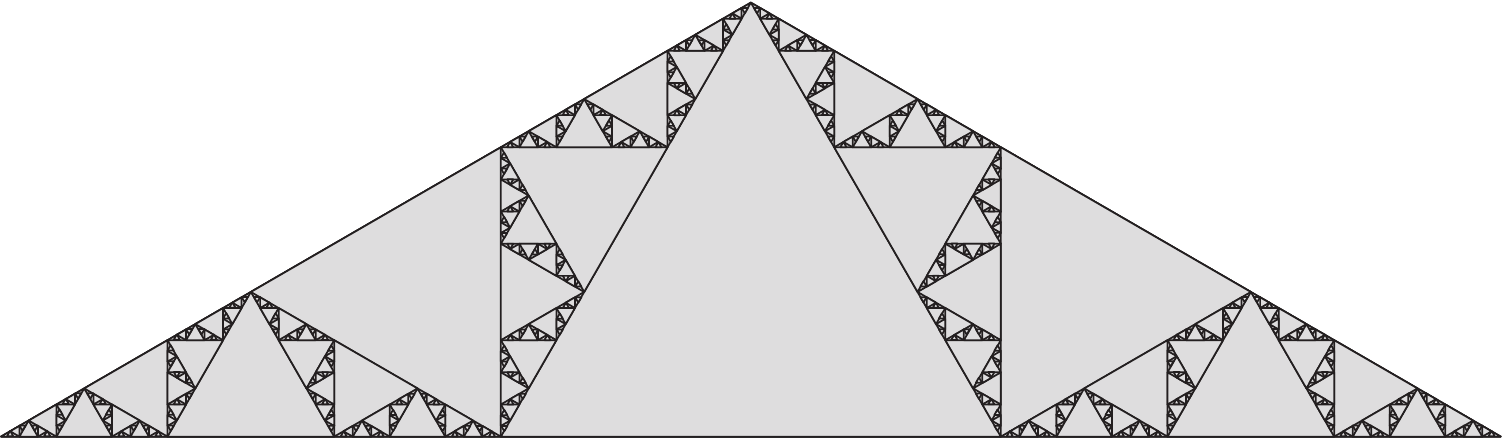}}
  \caption[The Koch tiling \sK.]
    {\captionsize The Koch tiling \sK.}
  \label{fig:koch-tiled-in}
  \centering
\end{figure}

The standard Koch tiling \sK (see Figure~\ref{fig:koch-tiled-in},
along with Figure~\ref{fig:koch-contractions} of \S\ref{sec:The
self-similar tiling}) is constructed via the self-similar system
\begin{align}
  \simt_1(z) := \gx \cj{z} \q \text{and} \q
  \simt_2(z) := (1-\gx)(\cj{z}-1)+1,
  \label{eqn:Koch system eqns}
\end{align}
with \(\gx = \tfrac{1}{2} + \tfrac{1}{2\sqrt{3}}\ii\) and $z \in
\bC$. The attractor of $\{\simt_1,\simt_2\}$ is the classical von
Koch curve. Thus \sK has one scaling ratio \(r=|\gx| = 1/\sqrt3,\)
and one generator $G$: an equilateral triangle of side length
$\tfrac13$ and generating inradius \(g=\tfrac{\sqrt3}{18}\). This
tiling has inradii $\gr_m = gr^{m}$ with multiplicity $2^m$, where $m=0,1,2,\dots$, so the
scaling zeta function is
\begin{align}
  \gzs(s) = \frac1{1 - 2 \cdot 3^{-s/2}},
\end{align}
and the scaling complex dimensions are
\begin{align}
  \Ds = \{D + \ii n\per \suth n \in \bZ\}
  \qq \text{for } D=\log_34, \; \per=\tfrac{4\gp}{\log3}.
\end{align}
By inspection, a tile with inradius $1/x$ will have tube formula
\begin{align}
  \label{eqn:Koch gtf}
  \gtf(x,\ge) =
    \begin{cases}
      3^{3/2} \left(-\ge^2 + 2\ge x\right), &\ge \leq 1/x, \\
      3^{3/2} x^2, &\ge \geq 1/x.
    \end{cases}
\end{align}
For fixed $x$, \eqref{eqn:Koch gtf} is clearly continuous at $\ge = 0^+$. Thus we have
\begin{align*}
  \gzs(s) &= \frac{1}{1 - 2 \cdot 3^{-s/2}} \qq\text{and} \\
  \crv_{0} &= -3^{3/2}, \crv_{1} = 2\cdot 3^{3/2}, \crv_{2} = -3^{3/2}.
\end{align*}
Now applying \eqref{eqn:single gen result}, the tube formula for
the Koch tiling \sK is
\begin{align}
  V({\sK},\ge)
  &= 3^{3/2} \genir^2 \sum_{\gw \in \Ds} \res{\frac{1}{1 - 2 \cdot 3^{-s/2}}}
    \left(-\tfrac1\gw + \tfrac{2}{\gw-1} - \tfrac{1}{\gw-2}\right)\left(\tfrac{\ge}{\genir}\right)^{2-\gw} \notag \\
    &\hstr[12]+ \tfrac\genir2 \gzs(0) \res[0]{-\tfrac{1}{s}} \left(\tfrac\ge\genir\right)^{2-0}
              + \tfrac\genir2 \gzs(1) \res[1]{\tfrac2{s-1}}  \left(\tfrac\ge\genir\right)^{2-1}
    \notag \\ 
  &= \frac{\genir}{\log3} \sum_{n \in \bZ}
    \left(-\tfrac1{D+\ii n\per} + \tfrac{2}{D-1+\ii n\per} - \tfrac{1}{D-2+\ii n\per}\right)\left(\tfrac{\ge}{\genir}\right)^{2-D-\ii n\per} \notag \\
    &\hstr[12]+ 3^{3/2} \ge^2 + \tfrac1{1-2\cdot3^{-1/2}} \ge,
  \label{eqn:Koch final in eps/genir}
\end{align}
where $D=\log_34$, \(g=\tfrac{\sqrt3}{18}\) and
$\per=\tfrac{4\gp}{\log3}$ as before.

\begin{remark}
  \label{rem:comparison of Koch tiling tube to Koch tube}
  In \cite{KTF}, a tube formula was obtained for the Koch curve itself (rather than for the tiling associated with it) and the possible complex dimensions of this curve were inferred to be
  \[\D_{\sK\star} = \{D+\ii n\per \suth n \in \bZ\} \cup \{0+\ii n\per \suth n \in \bZ\},\]
  where $D=\log_34$ and $\per=\frac{2\gp}{\log3}$. The line of poles above
  $D$ was expected\footnote{This set of complex dimensions was predicted in
  \cite{FGNT1}, \S 10.3.}, and agrees precisely with the results of
  this paper. The meaning of the line of poles above $0$ is still unclear.
  A zeta function for the Koch curve was not defined prior to the present
  paper; all previous reasoning was by analogy with \eqref{eqn:1-dim tube formula}.
\end{remark}

\begin{remark}[Nonlattice Koch tilings]
  \label{rem:nonlattice Koch}
  By replacing \(\gx = \tfrac{1}{2} + \tfrac{1}{2\sqrt{3}}\ii\)
  in \eqref{eqn:Koch system eqns} with any other complex number
  satisfying
  \(|\gx|^2 + |1-\gx|^2 < 1,\)
  one obtains a family of examples of nonlattice self-similar
  tilings. The tube formula computations parallel the lattice case
  almost identically. The lattice Koch tilings
  correspond exactly to those $\gx \in B(\tfrac12,\tfrac12)$ (the ball of
  radius $\tfrac12$ centered at $\tfrac12 \in \bC$) for which
  $\log_r |\gx|$ and $\log_r |1-\gx|$ are both positive integers,
  for some fixed $0<r<1$. Further discussion (and illustrations)
  of nonlattice Koch tilings may be found in \cite{SST}.
\end{remark}

\subsection{The Sierpinski gasket tiling}
\label{ssec:Sierpinski gasket tiling}

\begin{figure}
  \centering
  \scalebox{0.95}{\includegraphics{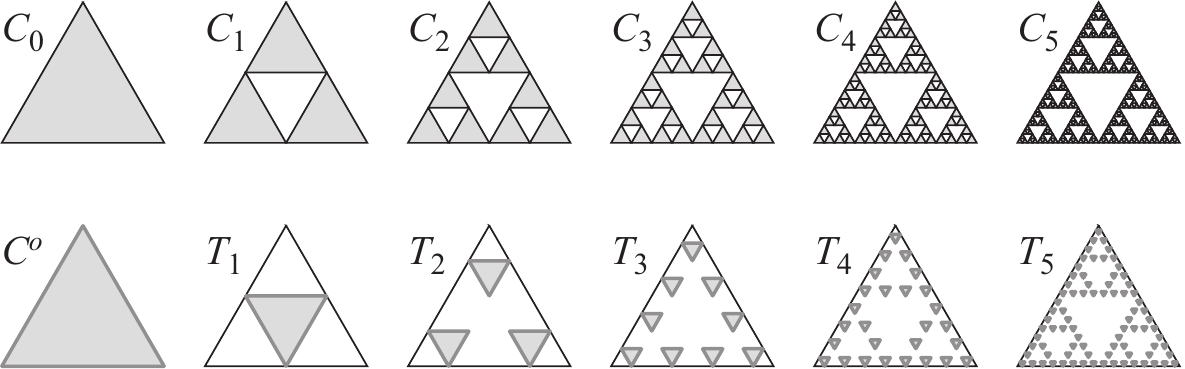}}
  \caption{The Sierpinski gasket tiling.}
  \label{fig:Sierpinski gasket tiling}
\end{figure}

The Sierpinski gasket tiling \sSG (see Figure~\ref{fig:Sierpinski
gasket tiling}) is constructed via the system
  \[\simt_1(z) := \tfrac12 z, \q
    \simt_2(z) := \tfrac12 z + \tfrac12, \q
    \simt_3(z) := \tfrac12 z + \tfrac{1+\ii\sqrt3}4,\]
which has one common scaling ratio \(r = 1/2,\) and one generator
$G$: an equilateral triangle of side length $\tfrac12$ and inradius
\(\genir = \tfrac1{4\sqrt3}.\) Thus \sSG has inradii $\gr_m = \genir
r^{m}$ with multiplicity $3^m$, $m=0,1,2,\dots$, so the scaling zeta function is
\begin{align}
  \gzs(s) = \frac1{1 - 3 \cdot 2^{-s}},
\end{align}
and the scaling complex dimensions are
\begin{align}
  \Ds = \{D + \ii n\per \suth n \in \bZ\}
  \qq \text{for } D=\log_23, \; \per=\tfrac{2\gp}{\log2}.
\end{align}

Aside from $\gzs(s)$, the tube formula calculation for \sSG is
identical to that for the previous example \sK:
\begin{align}
  V({\sSG},\ge)
  &= \tfrac{\sqrt3}{16 \log2} \sum_{n \in \bZ}
    \left(-\tfrac1{D+\ii n\per} + \tfrac{2}{D-1+\ii n\per} - \tfrac{1}{D-2+\ii n\per}\right)\left(\tfrac{\ge}{\genir}\right)^{2-D-\ii n\per} \notag \\
    &\hstr[12]+ \tfrac{3^{3/2}}{2} \ge^2 - 3 \ge.
  \label{eqn:Sierpinski final in eps/genir}
\end{align}

\begin{remark}
  \label{rem:epsnbd-for-noncolliders}
  Suppose that for a tiling \tiling, the boundary of the hull
  intersects the boundary of a generator in at most a finite set:
  \(  \left|\del \hull \cap \del \gen_q\right| < \iy.\)
  In this case, the tube formula for the tiling is almost the
  (exterior) tube formula for the attractor. This is the case for the
  Sierpinski gasket, and also for the Sierpinksi carpet (in which case
  the intersection 
  is empty). In
  fact, the exterior \ge-neighbourhood of the Sierpinski gasket curve is
  obtained by adding the Steiner's formula for $\cj{\hull}$:
  \begin{align}
    \vol[2]((\sSG)_\ge) = V({\sSG},\ge) + 3\ge + \gp\ge^2.
  \end{align}
\end{remark}

\subsection{The Pentagasket tiling} 
\label{exm:Pentagasket Tiling}

\begin{figure}
  \centering
  \scalebox{1.10}{\includegraphics{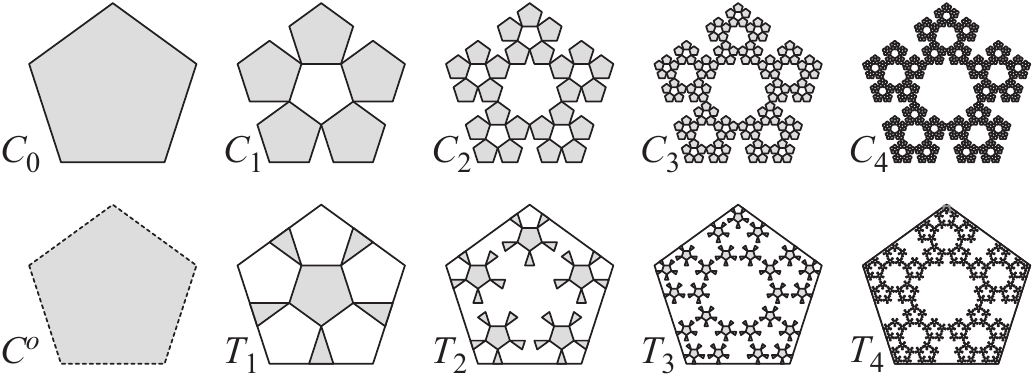}}
  \caption{The pentagasket tiling.}
  \label{fig:pentagasket tiling}
\end{figure}

The pentagasket tiling \sP (see Figure~\ref{fig:pentagasket tiling})
is constructed via the self-similar system defined by the five maps
\[\simt_\j(x) = \tfrac{3-\sqrt5}2 x + p_\j, \qq \j=1,\dots,5,\]
with common scaling ratio $r=\phi^{-2}$, where $\phi = (1+\sqrt5)/2$
is the golden ratio, and the points $\frac{p_\j}{1-r} = c_\j$ form
the vertices of a regular pentagon of side length 1.
%
\begin{figure}
  \centering
  \scalebox{0.85}{\includegraphics{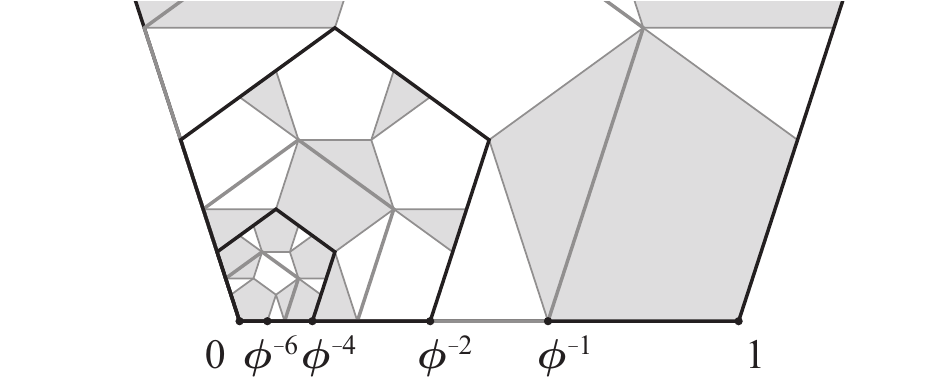}}
  \caption[The pentagasket and the golden ratio $\phi$.]
    {\captionsize The pentagasket and the golden ratio $\phi$.}
  \label{fig:pentagasket and golden ratio}
\end{figure}
%
The pentagasket \sP is an example of multiple generators
$\gen_q$: $\gen_1$ is a regular pentagon and
$\gen_2,\dots,\gen_6$ are congruent isosceles triangles, as
seen in $\tileset_1$ of Figure~\ref{fig:pentagasket tiling}. To
make the notation more meaningful, we use the subscripts $p,t$
to indicate a pentagon or triangle, respectively. The
generating inradius for the pentagon is $\genir_p =
\tfrac{\phi^2}2 \tan \tfrac3{10}\gp$ and the generating
inradius for the triangles is $\genir_t = \tfrac{\phi^3}2 \tan
\tfrac\gp5$. Thus, \sP has inradii $\gr_m = \genir_q r^{m}$,
for $q=p,t$ and $m=0,1,2,\dots$, with multiplicity $5^m$, so the scaling zeta function
is
\begin{align}
  \gzs(s) = \frac1{1 - 5 \cdot r^{-s}},
\end{align}
and the scaling complex dimensions are
\begin{align}
  \Ds = \{D + \ii n\per \suth n \in \bZ\}
  \qq \text{for } D=\log_{1/r}5, \; \per=\tfrac{2\gp}{\log r^{-1}}.
\end{align}
We omit the exercise of finding volumes for the pentagonal and
triangular generators; the tube formula for a tile of inradius $1/x$
is
\begin{align*}
  \gtfq(x,\ge) =
  \begin{cases}
    \crv_{q0}(\ge) x^0 + \crv_{q1}(\ge) x^1 =
    \ga_q \left(-\ge^2 + 2\ge x\right), &\ge \leq 1/x, \\
    \crv_{q2}(\ge) x^2 = \ga_q x^2, &\ge \geq 1/x,
  \end{cases}
\end{align*}
where $\ga_p:=5\cot\tfrac3{10}\gp$ and $\ga_t:=(\cot\tfrac\gp5)
/\left(1-\tan^2\tfrac\gp5\right)$. Since $\gen_2,\dots,\gen_6$
are congruent, we will apply Corollary~\ref{cor:self-similar-one-generator} to a triangle $\gen_t$ and multiply by 5 before
adding it to the result of applying
Corollary~\ref{cor:self-similar-one-generator} to the pentagon
$\gen_p$. For the pentagon and the triangle, we have $\crv_0 =
-\ga_q, \crv_1 = 2\ga_q$, and $\crv_2 = -\ga_q$.

The \tubezf zeta function of \sP is
\begin{align*}
  \gz_\sP(\ge,s)
  &= \sum_{q=1}^6 \frac{\ga_q\genir_q^s}{1 - 5 \cdot r^{-s}}
   \left(-\tfrac1s + \tfrac{2}{s-1} - \tfrac{1}{s-2}\right)\ge^{2-s},
\end{align*}
and the tube formula for \sP is
\begin{align}
  V({\sP},\ge)
  &= \frac{\ga_p}{\log r^{-1}} \sum_{n \in \bZ} \genir_p^2
    \left(-\tfrac1{D+\ii n\per} + \tfrac{2}{D-1+\ii n\per} - \tfrac{1}{D-2+\ii n\per}\right)\left(\tfrac{\ge}{\genir_p}\right)^{2-D-\ii n\per} \notag \\
    &\hstr[4]+ \frac{5\ga_t}{\log r^{-1}} \sum_{n \in \bZ} \genir_t^2
      \left(-\tfrac1{D+\ii n\per} + \tfrac{2}{D-1+\ii n\per} - \tfrac{1}{D-2+\ii n\per}\right)\left(\tfrac{\ge}{\genir_t}\right)^{2-D-\ii n\per} \notag \\
    &\hstr[4]+ \left[\left(\tfrac{\ga_p}4 + \tfrac{5\ga_t}4\right) \ge^2
      + \tfrac{(2\ga_p \genir_q + 10\ga_p \genir_q r)r}{r-5} \ge \right],
  \label{eqn:pentagasket final in eps/genir}
\end{align}
with $r=\phi^{-2}$, $\ga_p=5\cot\tfrac3{10}\gp$, $\ga_t=(\cot\tfrac\gp5) /\left(1-\tan^2\tfrac\gp5\right)$, $\genir_p = \tfrac{\phi^2}2 \tan \tfrac3{10}\gp$, $\genir_t = \tfrac{\phi^3}2 \tan \tfrac\gp5$, $D=\log_{1/r}5$ and \mbox{$\per=\tfrac{2\gp}{\log r^{-1}}$}.

\begin{remark}
  Much as in the case of fractal strings where $d=1$ (see \cite{FGCD}, \S8.4.2), it follows from Theorem~\ref{thm:self-similar case} that for a lattice self-similar tiling \tiling,
  each line of simple complex dimensions $\gb + \ii n\per$ gives rise to a function which consists of a multiplicatively periodic function times $\ge^{d-\gb}$. Here, \gb is some real constant and $\per = 2\gp/\log r^{-1}$ is the oscillatory period of \tiling. Consequently, since the scaling complex dimensions with real part $D$ are always simple, the tube formula for each lattice tiling in this section has the form
  \begin{align}
    V(\tiling,\ge)
    &= h\left(\log_{r^{-1}}(\ge^{-1})\right) \ge^{d-\gb} + P(\ge),
    \label{eqn:general multiplicative form}
  \end{align}
  where $h$ is an additively periodic function of period 1 and $P$ is a polynomial in \ge. For instance, the periodic function appearing in the tube formula \eqref{eqn:Koch final in eps/genir} for the Koch tiling \sK of Example~\ref{ssec:The Koch Tiling} has the following Fourier expansion:
  \begin{align}
    h(u)
    &= \tfrac{\genir}{\log 3} \sum_{n \in \bZ} \genir^{\ii n\per}
      \left(-\tfrac{1}{D+\ii n\per} + \tfrac2{D-1+\ii n\per} - \tfrac1{D-2+\ii n\per}\right)e^{2\gp \ii nu},
    \label{eqn:h(u) fourier expn for Koch}
  \end{align}
  where $\genir=\sqrt3/18$, $D=\log_34$, $r=1/\sqrt3$, and
  $\per=4\gp/\log3$.
  We note that multiplicatively periodic terms appear frequently in the mathematics and physics literature. See, for example, the relevant references given in \S1.5, \S2.7, \S6.6, and \S12.5 of \cite{FGCD}.
\end{remark}


\section{Some remarks on the results in this paper}
\label{sec:Concluding Remarks}

In this last section, we briefly comments on several consequences and possible extensions of our main results.




\begin{remark}
  The monograph \cite{FGCD} proposes a new definition of a
  \emph{fractal} as ``an object with nonreal complex dimensions that
  have a positive real part''.
  With respect to this definition, the present work confirms the fractal
  nature of all the examples discussed in \S\ref{sec:Tube formula examples},
  and more generally, of all self-similar tilings considered in this
  paper.
\end{remark}

\begin{remark}
  \label{rem:geom-interp-of-zeta(0)}
  Our results for tilings shed new light on the (\mbox{1-dimensional}) tube formula for fractal strings \eqref{eqn:form of 1-dim tube formula}. The origin of the previously mysterious linear term $\{2\ge\gzL(0)\}$ (see \eqref{eqn:1-dim tube formula})) is now seen to come from a monophase formula for the unit interval, akin to \eqref{def:monophase}. This is discussed further in \S\ref{ssec:The tube formula for fractal strings}.
  In fact, all terms coming from the second sum of the extended distributional formula of Theorem~\ref{thm:ext dist explicit formula} are now understood to be related to the curvatures of the generator. This reveals a geometric interpretation and allows the two sums to be naturally combined, as seen in \eqref{eqn:res at j} and \eqref{eqn:dissection of 2ge gz(0)}.
\end{remark}

\begin{remark}
  Many classical fractal curves are attractors of more than one
  self-similar system. For example, the Koch curve discussed in
  \S\ref{ssec:The Koch Tiling} is also the attractor of a system of four
  simlarity transformations of $\bR^2$, each with scaling ratio $r=\tfrac13$. In this particular
  example, changes in the scaling zeta function produce a different
  set of complex dimensions. In fact, we obtain a subset of the
  original complex dimensions: $\{\log_34 + \ii n\per \suth n \in \bZ,
  \per=4\gp/\log3\}$. This has a natural geometric interpretation
  which is to be discussed in later work. In particular, it would
  be desirable to determine precisely which characteristics remain
  invariant between different tilings which are so related.
\end{remark}

\begin{remark}\label{rem:obtaining-the-ext-from-the-int}
  The two formulas \eqref{eqn:preview-of-simplified-tube-formula} and \eqref{eqn:steiner_conceptual} initially appear to be measuring very different things, but this is misleading. If one considers the example of the Sierpinski tiling (discussed in \S\ref{ssec:Sierpinski gasket tiling}), then it is immediately apparent that the exterior \ge-neighbourhood of the Sierpinski gasket is, in fact, equal to the union of the inner \ge-neighbourhood of the tiling and the exterior \ge-neighbourhood of the largest triangle. With $C_0$ as in Figure~\ref{fig:Sierpinski gasket tiling},
  \begin{align}
    \vol[2](\sSG_\ge^{ext}) = V(\sSG,\ge) + \vol[2](C_0).
    \label{eqn:extnbd_decomp}
  \end{align}
  However, things do not always work out so neatly, as the example of the Koch tiling shows; see \S\ref{ssec:The Koch Tiling}. In the forthcoming paper \cite{GeometryOfSST}, precise conditions are given for equality to hold as in \eqref{eqn:extnbd_decomp}. This allows one to use results of the present paper to compute explicit tube formulas for a large family of self-similar sets, including the Sierpinski gasket and carpet; see Remark~\ref{ssec:Sierpinski gasket tiling}.

  Despite the fact that the tube formula for a self-similar tiling may differ from the tube formula for the corresponding self-similar set, it still gives us valuable information about self-similar geometries (and their associated dynamical systems). Indeed, we can define the complex dimensions of a given self-similar set in \bRd to be those of the self-similar tiling canonically associated to it (as in \cite{SST}). This is motivated by focusing on the dynamics of the self-similar system, rather than looking directly at the set. For an example, see \S\ref{ssec:The Koch Tiling}, especially Remark~\ref{rem:comparison of Koch tiling tube to Koch tube}. 
\end{remark}

\begin{remark}\label{rem:Weyl's-celebrated-formula}
  Recall that Weyl's celebrated formula \cite{We2} expresses the leading asymptotics of the eigenvalue counting function of the (Dirichlet) Laplacian on a $d$-dimensional compact Riemann manifold $M$ in terms of the volume of $M$ and its dimension $d$; see, e.g., \cite[\S12.5 and App.~B]{FGCD}. 
  An analogue of this formula exists for certain manifolds with fractal boundaries (for example, very irregular bounded open sets in \bRd), and in that case, the corresponding error term can be expressed in terms of the Minkowski dimension of the boundary; see \cite{La1,La2,La3}, along with \cite[\S12.5.1 and \S12.5.2]{FGCD} and the relevant references therein. For a fractal spray satisfying suitable hypotheses, we should be able to use the framework and the results of the present paper in order to obtain a full spectral asymptotic expansion (not just the leading term) for the Laplacian, expressed in terms of the underlying visible complex dimensions. 
   See \cite[\S6.3--\S6.5]{FGCD} for the 1-dimensional case, and \cite[\S.6.6]{FGCD} for an example in the case of a of self-similar fractal spray: the Sierpinski drum. 
   We hope to elaborate on this remark in a later work.
\end{remark}

\begin{remark}[Note added in proof]
  \label{rem:note-added-in-proof}
  It turns out that a little more care is needed for the results concerning Minkowski measurability. To ensure that $\lim_{\ge \to 0^+} V(\tiling,\ge) \ge^{-(d-D)}$ (or its counterpart in the lattice case, the average Minkowski content) does not vanish in \eqref{def:Mink meas}, and that the corresponding term in the explicit formula dominates the remaining ones.
  This issue was discovered during collaboration with Steffen Winter on the forthcoming paper \cite{Pointwise}, and so it will be discussed in full detail in that paper. For example, a sufficient condition for Corollary~\ref{cor:meas and dichotomy} to hold is that $D > d-1$ and 
  \begin{align}\label{Xeqn:M_D(G)}
    \sum_{k=0}^{d-1} \frac{\genir^{\abscissa-k}}{\abscissa-k} (d-k) \crv_k(\gen) \neq 0,
  \end{align}
  where \genir is the inradius of the monophase generator \gen and $\crv_k(\gen)$ is as in Definition~\ref{def:monophase}.
\end{remark}

\begin{appendix}

\section{The Definition and Properties of \gzT}
  \label{app:Rigourous Definition of gzT}

In this appendix, we confirm some basic properties of the
\tubezf zeta function \gzT of a fractal spray \tiling. However, we first
require some facts about Mellin transformation. If $\testfn \in
\bD=C_c^\iy(0,\iy)$, it is elementary to check that for every
$s \in \bC$, the Mellin transform $\testfnt(s)$ is given by the
well-defined integral \eqref{def:Mellin transform} and satisfies
$|\testfnt(s)| \leq |\testfnt|(\Re s) < \iy$. We will need
additional estimates in what follows. We also use the
forthcoming fact that $\testfnt(s)$ is an entire function; see Corollary~\ref{cor:holomorphicity of gft}.

Throughout Appendix~\ref{app:Rigourous Definition of gzT}, we assume that the hypotheses of Theorem~\ref{thm:fractal-spray-tube-formula} are satisfied.

\begin{lemma}
  \label{lem:uniform bound on gft}
  Suppose that $S \ci \bC$ is horizontally bounded, so $\inf S:= \inf_S \Re s$ and
  $\sup S := \sup_S \Re s$ are finite. Let $K$ be a compact interval
  containing the support of $\testfn \in C_c^\iy(0,\iy)$.
  Then there is a constant $c_K >0$ depending only on $K$ such that
  \begin{align}
    \sup_{s \in S} |\testfnt(s)| \leq c_K \|\testfn\|_\iy.
    \label{eqn:bound on gft_k}
  \end{align}
  In particular, $\testfnt(s)$ is always uniformly bounded on any screen $S$
  as in Definition~\ref{def:screen}.
\end{lemma}
  \begin{proof}
    Let $K$ be a compact interval containing the support of \testfn. Since
    \begin{align}
      |x^{s-1}| = x^{\Re s-1}
      \leq \begin{cases}
        x^{\sup S-1}, &x \geq 1, \\
        x^{\inf S-1}, &0 < x < 1,
      \end{cases}
      \label{eqn:bound on x^s}
    \end{align}
    one can define a bound
    \begin{align*}
      b_K := \sup_{x \in K} \max \{x^{\sup S-1},x^{\inf S-1}\}.
    \end{align*}
    Note that $b_K$ is finite because the function $x \mapsto \max\{x^{\sup S-1},
    x^{\inf S-1}\}$ is continuous on the compact set $K$, and hence is bounded.
    Then we use \eqref{eqn:bound on x^s} to bound \testfnt as follows:
    \begin{align}
      |\testfnt(s)|
      &\leq \int_0^\iy |x^{s-1}|\cdot|\testfn(x)| \, dx \notag \\
      &= \int_K x^{\Re s-1} |\testfn|(x) \, dx = \widetilde{|\testfn|}(\Re s)
        \label{eqn:derivation of bound on gft_k} \\
      &\leq b_K \|\testfn\|_\iy \cdot \vol(K). \notag
      \qedhere
    \end{align}
  \end{proof}

\begin{remark}
  \label{rem:estimate on hgft translate}
  The exact counterpart of Lemma~\ref{lem:uniform bound on gft}
  holds if $\testfnt(s)$ 
  is replaced by a translate $\testfnt(s-s_0)$, 
  for any $s_0 \in \bC$. Therefore,
  under the same assumptions as in Lemma~\ref{lem:uniform bound on gft},
  we have
  \begin{align}
    \label{eqn:bound on gft translate}
    \sup_{s \in S} |\testfnt(s-s_0)|
    \leq c_{K,s_0} \|\testfn\|_\iy,
  \end{align}
  where $c_{K,s_0} := b_{K,s_0} \cdot \vol(K)$, and
  \begin{align}
    \label{eqn:bound constant for gft translate}
    b_{K,s_0} := \sup_{x \in K} \max\{x^{\sup S - \Re s_0 -1},x^{\inf S - \Re s_0 -1}\} < \iy.
  \end{align}
  In particular, for any compact interval $K$ containing the support
  of \testfn, and for each fixed integer $k \geq 0$,
  \begin{align}
    \label{eqn:bound on_compact gft translate}
    \sup_{s \in S} \left|\testfnt(s-d+k+1)\right|
    \leq c_{K,k} \|\testfn\|_\iy,
  \end{align}
  where $c_{K,k}$ is a finite and positive constant.
\end{remark}

\begin{lemma}
  \label{lem:holomorphicity of F(s)}
  Let $(X,\gm)$ be a measure space.
  Define an integral transform by $F(s) = \int_X f(x,s) \,d\gm(x)$
  where
  \[|f(x,s)| \leq G(x), \q \text{for some $G \in L^1(X,\gm)$,}\]
  for \gm-a.e. $x \in X$, and for all $s$ in some neighbourhood of
  $s_0 \in \bC$.
  If the function \(s \mapsto f(x,s)\) is holomorphic for \gm-a.e.
  $x \in X$, then $F(s)$ is well-defined and holomorphic at $s_0$.
\end{lemma}

The proof is a well-known application of Lebesgue's Dominated
Convergence Theorem. We use Lemma~\ref{lem:holomorphicity of F(s)}
to obtain the following
corollary, which is used to prove Theorem~\ref{thm:fractal-spray-tube-formula} and Theorem~\ref{thm:gzT is a d-valued mero}.

\begin{cor}
  \label{cor:holomorphicity of gft}
  For $\testfn \in C_c^\iy(0,\iy)$, $\testfnt(s)$ is entire.
\end{cor}
  \begin{proof}
    Fix $s_0 \in \bC$. If $s$ is in a compact neighbourhood of $s_0$, then
    $\Re s$ is bounded, say by $\ga \in \bR$. Then for almost every
    $x>0$,
    \begin{align}
      \left|x^{s-1} \testfn(x)\right|
      \leq x^{\ga-1} \|\testfn\|_\iy \gc_\testfn,
      \label{eqn:}
    \end{align}
    where $\gc_\testfn$ is the characteristic function of the compact
    support of \testfn. Upon application of Lemma~\ref{lem:holomorphicity of F(s)},
    one deduces that \testfn is holomorphic at $s_0$.
  \end{proof}

Caution: Corollary~\ref{cor:holomorphicity of gft} does not combine with Lemma~\ref{lem:uniform bound on gft} to imply that \testfnt is constant; indeed, Liouville's Theorem does not apply here because $s$ is restricted to the screen $S$ in Lemma~\ref{lem:uniform bound on gft}.

\begin{defn}
  \label{def:d-valued mero}
  For $T(\ge,s)$ to be a \emph{weakly meromorphic distribution-valued function}
  on $W$, there must exist (i) a discrete set $\Poles \ci \bC$, and
  (ii) for each $\gw \in \Poles$, an integer $n_\gw < \iy$,
  such that $\Psi(s) = \la T(\ge,s),\testfn(\ge)\ra$ is a meromorphic function of $s \in W$,
  and each pole \gw of $\Psi$ lies in \Poles and has multiplicity
  at most $n_\gw$.

  To say that the distribution-valued function $T:W \to \Dist$ given by
  $s \mapsto T(\ge,s)$ is \emph{(strongly) meromorphic }means that, as a
  $\Dist$-valued function, it is truly a meromorphic function, in the sense
  of the proof of Lemma~\ref{lemma:weakly mero implies mero}.
  Recall that we are working with the space of distributions \Dist, defined
  as the dual of the space of test functions $\bD = C_c^\iy(0,\iy)$.
\end{defn}

\begin{lemma}
  \label{lemma:weakly mero implies mero}
  If $T$ is a weakly meromorphic distribution-valued function, then it is a (strongly) meromorphic
  distribution-valued function.
\end{lemma}
  \begin{proof}
    For $\gw \notin \Poles$, note that as $s \to \gw$,
    \begin{align}
      \label{eqn:diff quot of Ts}
      \frac{T(\ge,s)-T(\ge,\gw)}{s-\gw}
    \end{align}
    converges to a distribution (call it $T'(\ge,\gw)$) in \Dist, by the
    Uniform Boundedness Principle for a topological vector space
    such as \bD; see \cite{Rud3}, Thm.~2.5 and Thm.~2.8. Hence,
    the \Dist-valued function $T$ is holomorphic at \gw.

    For $\gw \in \Poles$, apply the same argument to
    \begin{align}
      \label{eqn:res limit of Ts}
      \lim_{s \to \gw} \frac1{(n_\gw-1)!} \left(\frac{d}{ds}\right)^{\negsp[1]n_\gw-1}
      \negsp[10]\left( (s-\gw)^{n_\gw} \vstr[2.2]T(\ge,s)\right),
    \end{align}
    which must therefore define a distribution, i.e., exist as an element
    of \Dist.
    Thus $T$ is truly a meromorphic function with values in \Dist, and with
    poles contained in \Poles.
  \end{proof}

\begin{theorem}
  \label{thm:gzT is a d-valued mero}
  Under the hypotheses of Theorem~\ref{thm:fractal-spray-tube-formula}
  or Theorem~\ref{thm:self-similar case}, the \tubezf zeta function
  of a fractal spray or tiling
  \begin{align}
    \label{eqn:gzT as double sum}
    \gzT(\ge,s)
    &= \ge^{d-s} \gzs(s) \sum_{k=0}^d  \frac{\genir^{s-k}}{s-k} \crv_{k}
  \end{align}
  is a distribution-valued (strongly) meromorphic function on $W$, with poles
  contained in \DT.
\end{theorem}
  \begin{proof}
    Let $\Poles=\DT$ and note that
    \begin{align}
      \label{eqn:gzT is mero for any gf}
      \la \gzT(\ge,s),\testfn(\ge)\ra
      &= \gzs(s) \sum_{k=0}^d \frac{\genir^{s-k}}{s-k} \crv_k \int_0^\iy \ge^{d-s} \testfn(\ge) \, d\ge \notag \\
      &= \gzs(s) \testfnt(d-s+1) \sum_{k=0}^d \frac{\genir^{s-k}}{s-k} \crv_k.
    \end{align}
    By Corollary~\ref{cor:holomorphicity of gft}, this is a finite sum of
    meromorphic functions and hence meromorphic on $W$, for any test function \testfn.
    Applying Lemma~\ref{lemma:weakly mero implies mero}, one sees
    that \gzT is a meromorphic function with values in \Dist.
%
  \end{proof}

\begin{remark}
  \label{rem:bounds on the multiplicity of poles}
  Note that for each $\testfn \in \bD$, the poles of the
  \bC-valued function
  \begin{align}
    \label{eqn:bounds on the multiplicity of poles}
    s \mapsto \left\la \gzT(\ge,s),\testfn(\ge)\right\ra
  \end{align}
  are contained in \DT. Further, if $m_\gw$ is the multiplicity
  of $\gw \in \DT$ as a pole of $\gzs(s)$, then the multiplicity
  of \gw as a pole of \eqref{eqn:bounds on the multiplicity of poles}
  is bounded by $m_\gw +1$.
\end{remark}

\begin{cor}
  \label{cor:res(gzT) is a wdef dist}
  The residue of \gzT at a pole $\gw \in \DT$ is a well-defined distribution.
\end{cor}
  \begin{proof}
    This follows immediately from the second part of the proof of
    Lemma~\ref{lemma:weakly mero implies mero}, with $\Poles=\DT$.
%
%
  \end{proof}

\begin{cor}
  The sum of residues appearing in  Theorem~\ref{thm:fractal-spray-tube-formula} and Theorem~\ref{thm:self-similar case} is distribution{ally} convergent,  and is thus a well-defined distribution.
\end{cor}
  \begin{proof}
    In view of the proof of Theorem~\ref{thm:gzT is a d-valued mero},
    this comes by applying the Uniform Boundedness Principle to
    an appropriate sequence of partial sums, in a manner similar
    to the proof of Lemma~\ref{lemma:weakly mero implies mero}. Again,
    see \cite[Rem.~5.21]{FGCD}.
  \end{proof}


\section{The Error Term and Its Estimate}
\label{app:The Error Term and Estimate}

In this appendix, we give the promised proof of the expression for
the error term \eqref{def:main error term} and its estimate \eqref{eqn:main error est}, as stated in Theorem~\ref{thm:fractal-spray-tube-formula}. 
Throughout Appendix~\ref{app:The Error Term and Estimate}, we assume that the hypotheses of Theorem~\ref{thm:fractal-spray-tube-formula} are satisfied. First, we require a definition.

\begin{defn}[Primitives of distributions]
  \label{def:primitive of a dist}
  Let $T_\gh$ be a distribution defined by a measure as
    $\la T_\gh, \testfn \ra := \int_0^\iy \testfn \,d\gh.$
  Then the \kth \emph{primitive} (or \kth antiderivative) of $T_\gh$ is
  defined by $\la T_\gh^{[k]},\testfn\ra := (-1)^k \la T_\gh,\testfn^{[k]}
  \ra,$
  where $\testfn^{[k]}$ is the \kth primitive of $\testfn \in C_c^\iy(0,\iy)$
  that vanishes at \iy \xspace together with all its derivatives.
  For $k \geq 1$, for example,
  \begin{align}
    \label{eqn:kth primitive of T as an integral}
    \la T_\gh^{[k]},\testfn\ra = \int_0^\iy \int_y^\iy \frac{(x-y)^{k-1}}{(k-1)!} \testfn(x) \, dx \, d\gh(y).
  \end{align}
\end{defn}

\begin{theorem}
  \label{thm:Mellin transform of the kth primitive}
  The Mellin transform of the \kth primitive of a test function is
  given by $\widetilde{\testfn^{[k]}}(s) = \tilde\testfn(s+k) \gx_k(s)$,
  where $\gx_k$ is the meromorphic function
  \begin{align}
    \gx_k(s) := \sum_{j=0}^{k-1} \frac{\binom{k-1}j \raisebox{-0.9mm}{$(-1)^j$}}{(k-1)! (s+j)}.
    \label{eqn:defn gy_k}
  \end{align}
\end{theorem}
  \begin{proof}
    By direct computation,
    \begin{align}
      \widetilde{\testfn^{[k]}}(s)
      &= \int_0^\iy \ge^{s-1} \int_\ge^\iy \frac{(x-\ge)^{k-1}}{(k-1)!} \testfn(x) \,dx \,d\ge \notag \\
      &= \frac1{(k-1)!} \int_0^\iy \int_\ge^\iy \sum_{j=0}^{k-1} \binom{k-1}j x^{k-1-j} (-\ge)^j \ge^{s-1} \testfn(x) \,dx \,d\ge \notag \\
      &= \sum_{j=0}^{k-1} \frac{\binom{k-1}j \raisebox{-0.9mm}{$(-1)^j$}}{(k-1)!} \int_0^\iy \int_\ge^\iy x^{k-1-j} \ge^{s+j-1} \testfn(x) \,dx \,d\ge \notag \\
      &= \sum_{j=0}^{k-1} \frac{\binom{k-1}j \raisebox{-0.9mm}{$(-1)^j$}}{(k-1)!} \int_0^\iy x^{k-1-j} \testfn(x) \int_0^x \ge^{s+j-1} \,d\ge \,dx \notag \\
      &= \sum_{j=0}^{k-1} \frac{\binom{k-1}j \raisebox{-0.9mm}{$(-1)^j$}}{(k-1)! (s+j)} \int_0^\iy x^{s+k-1} \testfn(x) \,dx
        \label{eqn:the birth of gy_k} \\
      &= \tilde\testfn(s+k) \gx_k(s). \notag
    \end{align}
    Again, the formula \eqref{eqn:defn gy_k} for $\gx_k$ is valid for
    $\Re s > k$ by \eqref{eqn:the birth of gy_k}, but then extends to
    being valid for all $s \in \bC$ by meromorphic continuation.
  \end{proof}

\begin{cor}
  \label{cor:Mellin transform of kth primitive}
  We also have $\left|\widetilde{\testfn^{[k]}}(s)\right|
  \leq \left|\testfnt(s+k) \gx_k(s)\right|$.
\end{cor}

\begin{remark}
  \label{rem:bound on gx_k}
  For $s \in \bC, t=\Im s$, and $c_{\gx}>0$, we also have
  \begin{align}
    |\gx_k(s)| \leq \frac{c_\gx}{|t|^k}.
    \label{eqn:bound on gx_k}
  \end{align}
\end{remark}

We are now in a position to provide the proofs previously
promised.
\begin{theorem}
  \label{thm:validation of error dist}
  As stated in \eqref{def:main error term} of
  Theorem~\ref{thm:fractal-spray-tube-formula}, the error term is given by
  \begin{align}
    \label{eqn:error term for proof}
    \sR(\ge) 
    = \frac1{2\gp \ii}\int_S \gzT(\ge,s) \, ds,
  \end{align}
  and is a well-defined distribution.
\end{theorem}
  \begin{proof}
    Applying \eqref{def:Mellin transform} to \eqref{def:partial error Rq}
    for $k=0,\dots,d$ gives\footnote{In the proof of
    Theorem~\ref{thm:fractal-spray-tube-formula},
    the quantity \eqref{eqn:Mellin first applied to error} was denoted
    by $\la\sR_{k}, \testfn\ra$, so
    that \sR could easily be written (formally) as a function
    in \eqref{eqn:V as a function}. For clarity, since we work with test
    functions, this quantity is instead denoted by $\la\sR,
    \testfn\ra_{k}$ throughout this proof.}
    \begin{align}
       \label{eqn:Mellin first applied to error}
       \left\la\sR, \gf\right\ra_{k}
        =& \frac{1}{2\gp \ii} \int_S \frac{\genir^{s-k}}{s-k} \gzs(s) \crv_{k}
        \int_0^\iy \ge^{d-s} \testfn(\ge) \,d\ge
        \,ds.
    \end{align}
    To see that this gives a well-defined distribution \sR, we apply the descent method, as described in \cite{FGCD}, Rem.~5.20. The first step is to show that $\left\la\sR^{[k]}, \testfn\right\ra_{k}$ is a well-defined distribution for sufficiently large $k$; specifically, for any integer $k>\languidOrder$, where \languidOrder is the order of languidity, as in Definition~\ref{def:languid}. Note that we can break the integral along the screen $S$ into two pieces and work with each separately:
    \begin{align}
       \left\la\sR^{[k]}, \testfn\right\ra_{k}
        &= \frac{(-1)^k }{2\gp \ii} \int_{|\Im s| > 1} \frac{\genir^{s-k}}{s-k} \gzs(s) \crv_{k} \int_0^\iy \ge^{d-s} \testfn^{[k]}(\ge) \,d\ge \,ds
         \label{eqn:error R(x) impart pos}
         \\
        &+ \frac{(-1)^k }{2\gp \ii} \int_{|\Im s| \leq 1} \frac{\genir^{s-k}}{s-k} \gzs(s)\crv_{k}  \int_0^\iy \ge^{d-s} \testfn^{[k]}(\ge) \,d\ge \,ds.
         \label{eqn:error R(x) impart mid}
    \end{align}
    Here and throughout the rest of this appendix, it is understood
    that such integrals (as in \eqref{eqn:error R(x) impart pos}--\eqref
    {eqn:error R(x) impart mid}) are for $s \in S$.
    Since the screen avoids the integers $0,\dots,d$ by assumption,
    the quantity $|s-k|$ is bounded away from 0.
    Since the screen avoids the poles of $\gzs$ by hypothesis,
    $\gzs(s)$ is continuous on the compact set $\{s \in S \suth |\Im s|
    \leq 1\}$. Therefore, it is clear that
    \eqref{eqn:error R(x) impart mid} is a well-defined integral. We focus
    now on \eqref{eqn:error R(x) impart pos}:
    \begin{align}
      \left| \frac{\crv_{k}}{2\gp \ii}\right. & \left. \int_{|\Im s| > 1} \frac{\genir^{s-k}}{s-k} \gzs(s) \int_0^\iy \ge^{d-s} \testfn^{[k]}(\ge) \,d\ge \,ds \right| \notag \\
       &\leq \frac{\crv_{k}}{2\gp} \int_{\Im s > 1} \left| \genir^{s-k} \frac{\gzs(s)}{s-k} \right| \cdot \left|\widetilde{\testfn^{[k]}}(s-d+1)\right| \,ds \notag \\
       &\leq c_1 \int_1^\iy |t|^{M-1}
        \cdot \left|\testfnt(s - d + k + 1)\right| \cdot \left|\gx_k(s-d+1)\right| \,dt \notag \\
       &\leq c_1 \int_1^\iy c_{k} |t|^{M-1} \cdot c_K \|\testfn\|_\iy \cdot \frac{c_\gx}{|t|^k} \,dt, \notag \\
       &= C \|\testfn\|_\iy \int_1^\iy |t|^{M-1-k} \,dt,
       \label{eqn:bound on sR^k}
    \end{align}
    which is clearly convergent for $k>M$. The second inequality in \eqref{eqn:bound on sR^k} comes by condition \textbf{L2} of Definition~\ref{def:languid}. Also, recall (from the remark just after the statement of Lemma~\ref{lem:uniform bound on gft}) that for $s \in S$, the real part of $s$ is given by some function $f$ which is Lipschitz (cf.~Definition~\ref{def:screen}), and hence is almost everywhere differentiable and has a bounded derivative (where it exists) on the support of \testfn. The third comes by inequality \eqref{eqn:bound constant for gft translate} of Remark~\ref{rem:estimate on hgft translate}, along with Remark~\ref{rem:bound on gx_k}. This establishes the validity of $\la\sR^{[k]}, \testfn\ra_{k}$ and thus shows that $\sR^{[k]}$ defines a linear functionl on \bD.

    To check that the action of $\sR^{[k]}$ is continuous on \bD,
    let $\testfn_n \to 0$ in \bD, so that there is a compact set $K$ which contains the
    support of every $\testfn_n$, and $\|\testfn_n\|_\iy \to 0$. Then
    \begin{align}
      \left| \la\sR^{[k]}, \testfn_n \ra \right|
      &\leq C \cdot \left|\testfnt_n(s - d + k + 1)\right|
      \leq c_K \|\testfn_n\|_\iy \limas{n} 0,
       \label{eqn:continuity of sR^k}
    \end{align}
    by following \eqref{eqn:bound on sR^k} and then applying
    Lemma~\ref{lem:uniform bound on gft}, along with its extensions
    as stated in Remark~\ref{rem:estimate on hgft translate}.
    Thus, $\sR^{[k]}$ is a well-defined distribution.
    If we differentiate it distributionally $k$ times, we obtain \sR. This shows
    that \sR is a well-defined distribution and concludes the proof.
  \end{proof}

Before finally checking the error estimate, we define what is
meant by the expression $T(x)=O(x^\ga)$ as $x \to \iy$, when $T$
is a distribution.

\begin{defn}
  \label{def:order of a dist}
  When $\sR(x)=O(x^\ga)$ as $x \to \iy$ (as in \eqref{def:error term estimate}), we say as in \cite{FGCD}, \S5.4.2, that \sR is of \emph{asymptotic order at most $x^\ga$} as $x \to \iy$. To understand this expression, first define
  \begin{align}
    \label{def:gf_a}
    \testfn_a(x) := \tfrac1a \testfn\left(\tfrac{x}a\right),
  \end{align}
  for $a>0$ and for any test function \testfn. Then ``$\sR(x)=O(x^\ga)$ as $x \to \iy$'' means that
  \[\la\sR,\testfn_a\ra=O(a^\ga), \qq \text{ as } a \to \iy,\]
  for every test function \testfn. The implied constant may depend on \testfn. Similarly, ``$\sR(x)=O(x^\ga)$ as $x \to 0^+$'' (as in \eqref{eqn:1-dim error bound} and \eqref{eqn:main error est}) is defined to mean that
  \[\la\sR,\testfn_a\ra=O(a^\ga), \qq \text{ as } a \to 0^+,\]
  for every test function \testfn.
\end{defn}

\begin{theorem}[Error estimate]
  \label{rem:Error term estimate}
  As stated in Theorem~\ref{thm:fractal-spray-tube-formula},
  the error term $\sR(\ge)$ in \eqref{eqn:error term for proof}
  is estimated by
  \begin{align}
    \label{eqn:main error est in app}
    \sR(\ge) &=  O(\ge^{d-\sup S}), \qq \text{as } \ge \to 0^+.
  \end{align}
\end{theorem}
  \begin{proof}
    As in the proof of Theorem~\ref{thm:validation of error dist}, we use the
    descent method and begin by splitting the integral into two pieces. Since
    \[\la \sR^{[k]}, \testfn_a \ra = (-1)^k \la \sR, (\testfn_a)^{[k]} \ra,\]
    we work with 
    \begin{align}
       \left\la\sR, (\testfn_a)^{[k]}\right\ra_{k}
        &= \frac{\crv_{k}}{2\gp \ii}
           \int_{|\Im s| > 1} \negsp[5] \frac{\genir^{s-k}}{s-k} \gzs(s)
           \int_0^\iy \ge^{d-s} (\testfn_a)^{[k]}(\ge) \,d\ge \,ds
         \label{eqn:error est impart pos}
         \\
        &+ \frac{\crv_{k}}{2\gp \ii}
           \int_{|\Im s| \leq 1} \negsp[5] \frac{\genir^{s-k}}{s-k} \gzs(s)
           \int_0^\iy \ge^{d-s} (\testfn_a)^{[k]}(\ge) \,d\ge \,ds.
         \label{eqn:error est impart mid}
    \end{align}
    The \kth primitive of $\testfn_a$ is given by
    \begin{align}
      (\testfn_a)^{[k]}(\ge)
      &= \int_\ge^\iy \frac{(u-\ge)^{k-1}}{(k-1)!} \tfrac1a\testfn\left(\tfrac ua\right)\,du 
      = \int_{\ge/a}^\iy \frac{(au-\ge)^{k-1}}{(k-1)!} \testfn(u)\,du.
      \label{eqn:kth primitive of gf_a}
    \end{align}
    By following the same calculations as in
    Theorem~\ref{thm:Mellin transform of the kth primitive},
    one observes that
    \begin{align}
      &\left|\int_0^\iy  \frac{\ge^{d-s}}{s-k} \int_{\ge/a}^\iy \frac{(au-\ge)^{k-1}}{(k-1)!} \testfn(u) \,du \,d\ge \right|\notag \\
      &= \left| \int_0^\iy \int_0^{au} \frac{}{s-k} \sum_{j=0}^{k-1} \frac{\binom{k-1}j\raisebox{-0.9mm}{$(-1)^j$}}{(k-1)!} (au)^{k-1-j} \ge^{d-s+j} \testfn(u) \, d\ge \,du \right| \notag \\
      &\leq \frac{1}{|s-\i|} \sum_{j=0}^{k-1} \frac{\binom{k-1}j \raisebox{-0.9mm}{$(-1)^j$}}{(k-1)!} \int_0^\iy \left| (au)^{k-1-j} \testfn(u) \right| \int_0^{au} \left|  \ge^{d-s+j} \, d\ge \right| \,du \notag \\
      &\leq \frac{c_{k}}{|s-\i|} \gx_k(d-\Re s+1) \int_0^\iy (au)^{k-1-j} (au)^{d-\Re s+j+1} |\testfn(u)| \,du \notag \\
      &= a^{d-\Re s+k} \frac{c_{k}}{|s-\i|}  \gx_k(d-\Re s+1) \widetilde{|\testfn|}(d-\Re s+k) .
      \label{eqn:bound for crvqi with a}
    \end{align}
    Using \eqref{eqn:bound on gx_k} for $\gx_k$ and
    \eqref{eqn:derivation of bound on gft_k} for $\widetilde{|\testfn|}$
    (see Remark~\ref{rem:estimate on hgft translate}),
    we bound \eqref{eqn:error est impart pos} by
    \begin{align}
      \frac{c_{k}}{2\gp} &\int_{|\Im s| > 1} \negsp a^{d-\Re s+k} \cdot \frac{|\genir^{s-k}  \gzs(s)|}{|s-i|} \cdot \frac{c_\gy}{|t|^k} \cdot c_K \|\testfn\|_\iy  \,ds \\
      & \leq a^{d-\sup S +k} \left(C \int_1^\iy |t|^{M-1-k} \,dt\right),
      \label{eqn:bound on pos err est}
    \end{align}
    for any $0<a<1$, as in \eqref{eqn:bound on sR^k}.
    Since the integral in
    \eqref{eqn:bound on pos err est} clearly converges for $k>M$,
    we have established the estimate for $\sR^{[k]}$, along the
    part of the integral where $|\Im s|>1$. Recall that all our
    contour integrals are taken along the screen $S$.
    The proof for \eqref{eqn:error est impart pos}, where
    $|\Im s|>1$, readily follows from the corresponding argument
    in the proof of Theorem~\ref{thm:validation of error dist}. Thus
    we have established that
    \begin{align}
      \left|\la \sR^{[k]}(\ge), \testfn_a(\ge) \ra\right| \leq a^{d-\sup S+k} c_k,
      \q \text{for all } 0<a<1.
      \label{eqn:bound on Rk acting on phi_a}
    \end{align}
    In \eqref{eqn:bound on Rk acting on phi_a}--\eqref{eqn:bound on R acting on phi_a},
    the constants $c_k$ may depend on the test function
    \testfn.\footnote{Note that $c_{k-1}$ does not correspond to $c_k$
    when $k$ is replaced by $k-1$; rather, $c_{k-1}$ depends on
    the support of $\testfn'$. The notation is just used to indicate
    the analogous roles the constants $c_k$ play.}

    By iterating the following calculation:
    \begin{align}
      \left|\la \sR^{[k-1]}(\ge), \testfn_a(\ge) \ra\right|
      &= \left|\left\la \sR^{[k]}(\ge), \left(\tfrac1a \testfn\left(\tfrac\ge a\right)\right)' \right\ra\right| \notag \\
      &= \left|\tfrac1a \la \sR^{[k]}(\ge), (\testfn')_a (\ge) \ra\right| \notag \\
      &\leq a^{d-\sup S+k-1} c_{k-1},
      \label{eqn:bound on Rk-1 acting on phi_a}
    \end{align}
    one sees that
    \begin{align}
      \left|\la \sR(\ge), \testfn_a(\ge) \ra\right| \leq a^{d-\sup S}c_0,
      \q \text{for all } 0<a<1.
      \label{eqn:bound on R acting on phi_a}
    \end{align}
    By Definition~\ref{def:order of a dist}, this implies that $\sR(\ge) =
    O(\ge^{d-\sup S})$ as $\ge \to 0^+$.
  \end{proof}

\end{appendix}
 

\pgap

\par


\begin{thebibliography}{widest-label}
     \pgap[0.5cm]

     \bibitem{BeGo}
     \bibbook{M. Berger and B. Gostiaux}
        {Differential Geometry: Manifolds, Curves and Surfaces}
        {English transl., Springer-Verlag}
        {Berlin}
        {1988}



     \pgap[\bibgap]

     \bibitem{Fal1}
     \bibbook{K. J. Falconer}
        {Fractal Geometry --- Mathematical Foundations and Applications}
        {John Wiley}
        {Chichester}
        {1990}

     \pgap[\bibgap]

     \bibitem{Fal2}
     \bibart{K. J. Falconer}
        {On the Minkowsi measurability of fractals}
        {Proc. Amer. Math. Soc.}
        {}{123}{1995}{1115--1124}

     \pgap[\bibgap]

     \bibitem{Fed}
     \bibart{H. Federer}
        {Curvature measures}
        {Trans. Amer. Math. Soc.}
        {}{93}{1959}{418--491}

     \pgap[\bibgap]

     \bibitem{Fra}
     \bibart{M. Frantz}
        {Minkowski measurability and lacunarity of self-similar sets in \bR}
        {Proc. Symposia Pure Math.,}
        {Part 1,}{72}{2004}{77--91}

     \pgap[\bibgap]

     \bibitem{Fu1}
     \bibart{J. H. G. Fu}
        {Tubular neighbourhoods in Euclidean spaces}
        {Duke Math. J.}
        {}{52}{1985}{1025--1046}

     \pgap[\bibgap]

     \bibitem{Fu2}
     \bibart{J. H. G. Fu}
        {Curvature measures of subanalytic sets}
        {Amer. J. Math.}
        {}{116}{1994}{819--880}

     \pgap[\bibgap]

     \bibitem{Gat}
     \bibart{D. Gatzouras}
        {Lacunarity of self-similar and stochastically self-similar sets}
        {Trans. Amer. Math. Soc.}
        {}{352}{2000}{1953--1983}

     \pgap[\bibgap]

     \bibitem{Gr}
     \bibbook{A. Gray}
        {Tubes \emph{(2nd ed.)}}
        {Progress in Math., vol. 221, Birkh\"{a}user}
        {Boston}
        {2004}

     \pgap[\bibgap]

     \bibitem{HaLa}
     \bibart{B. M. Hambly and M. L. Lapidus}
        {Random fractal strings: Their zeta functions, complex dimensions and spectral asymptotics}
        {Trans. Amer. Math. Soc.}
        {}{358}{2006}{285--314}

     \pgap[\bibgap]

     \bibitem{HeLa1}
     \bibart{C. Q. He and M. L. Lapidus}
       {Generalized Minkowski content and the vibrations of fractal drums
         and strings}
       {Mathematical Research Letters}
       {} {3} {1996}{31-–40}

     \pgap[\bibgap]

     \bibitem{HeLa2}
     \bibart{C. Q. He and M. L. Lapidus}
       {Generalized Minkowski content, spectrum of fractal drums,
         fractal strings and the Riemann zeta-function}
       {Memoirs Amer. Math. Soc.}
       {No. 608,}{127}{1997}{1–-97}

     \pgap[\bibgap]

     \bibitem{HuLaWe}
     \bibart{D. Hug, G. Last and W. Weil}
        {A local Steiner-type formula for general closed sets and applications}
        {Mathematische Zeitschrift.}
        {}{246}{2004}{237--272}

     \pgap[\bibgap]

    \bibitem{Hut}
    \bibart{J. E. Hutchinson}
        {Fractals and self-similarity}
        {Indiana Univ. Math. J.}
        {}{30}{1981}{713--747}

     \pgap[\bibgap]

    \bibitem{Kig}
    \bibbook{J. Kigami}
        {Analysis on Fractals}
        {Cambridge Univ. Press}
        {Cambridge}
        {1999}

     \pgap[\bibgap]

    \bibitem{Rota}
    \bibbook{D. A. Klain and G.-C. Rota}
        {Introduction to Geometric Probability}
        {Accademia Nazionale dei Lincei, Cambridge Univ. Press}
        {Cambridge}
        {1999}

     \pgap[\bibgap]

     \bibitem{La1}
     \bibart{M. L. Lapidus}
        {Fractal drum, inverse spectral problems for elliptic operators and a partial resolution of the Weyl--Berry conjecture}
        {Trans. Amer. Math. Soc.}
        {}{325}{1991}{465--529}

     \pgap[\bibgap]

     \bibitem{La2}
        M. L. Lapidus,
        \emph{Spectral and fractal geometry: From the Weyl--Berry conjecture
        for the vibrations of fractal drums to the Riemann zeta-function,}
        in: Differential Equations and Mathematical Physics
        (C. Bennewitz, ed.),
        Proc. Fourth UAB Internat. Conf. (Birmingham, March 1990),
        Academic Press, New York, 1992, pp. 151--182.

     \pgap[\bibgap]

     \bibitem{La3}
        M. L. Lapidus,
        \emph{Vibrations of fractal drums, the Riemann hypothesis, waves in fractal
        media, and the Weyl--Berry conjecture,}
        in: Ordinary and Partial Differential Equations
        (B. D. Sleeman and R. J. Jarvis, eds.),
        vol. IV, Proc. Twelfth Internat. Conf.
        (Dundee, Scotland, UK, June 1992),
        Pitman Research Notes in Math. Series, vol. 289,
        Longman Scientific and Technical, London, 1993, pp. 126--209.

     \pgap[\bibgap]

     \bibitem{LaMa}
     \bibart{M. L. Lapidus and H. Maier}
        {The Riemann hypothesis and inverse spectral problems for fractal strings}
        {J. London Math. Soc.}
        {(2)}{52}{1995}{15--34}

     \pgap[\bibgap]

     \bibitem{KTF}
     \bibart{M. L. Lapidus and E. P. J. Pearse}
        {A tube formula for the Koch snowflake curve, with applications to complex dimensions}
        {J. London Math. Soc.}
        {(2) No. 2,}{74}{2006}{397--414. \arxiv{math-ph/0412029}}

     \pgap[\bibgap]

     \bibitem{TFSST}
     {M. L. Lapidus and E. P. J. Pearse,}
        {Tube formulas for self-similar fractals, in \emph{Analysis on Graphs and Its Applications} (P. Exner, J. P. Keating, C. Bristol, P. Kuchment, T. Sunada, and A. Teplyaev, eds.), Proc. of Symposia in Pure Mathematics, vol. 77, Amer. Math. Soc., Providence, RI, 2008, pp. 211--230. \arxiv{0711.0173}.}






     \pgap[\bibgap]

     \bibitem{FCM}
     M. L. Lapidus, E. P. J. Pearse and S. Winter,
        {Fractal curvature measures and local tube formulas,}
        {work in progress.}

     \pgap[\bibgap]

     \bibitem{Pointwise}
     M. L. Lapidus, E. P. J. Pearse and S. Winter,
        {Pointwise tube formulas for fractal sprays and self-similar tilings with arbitrary generators,}
        {in preparation.}

     \pgap[\bibgap]

     \bibitem{LaPo1}
     \bibart{M. L. Lapidus and C. Pomerance}
        {The Riemann-zeta function and the one-dimensional Weyl--Berry conjecture for fractal drums}
        {Proc. London Math. Soc.}
        {(3)}{66}{1993}{41--69}

     \pgap[\bibgap]

     \bibitem{LaPo2}
     \bibart{M. L. Lapidus and C. Pomerance}
        {Counterexamples to the modifed Weyl--Berry conjecture on fractal drums}
        {Math. Proc. Cambridge Philos. Soc.}
        {}{119}{1996}{167--178}

     \pgap[\bibgap]

     \bibitem{FGNT1}
     \bibbook{M. L. Lapidus and M. van Frankenhuysen}
        {Fractal Geometry and Number Theory: Complex dimensions of fractal strings and zeros of zeta functions}
        {Birkh\"{a}user}
        {Boston}
        {2000}

     \pgap[\bibgap]

     \bibitem{La-vF3}
     \bibart{M. L. Lapidus and M. van Frankenhuysen}
        {Complex dimensions of self-similar fractal strings and Diophantine approximation}
        {J. Experimental Mathematics}
        {No. 1,}{12}{2003}{41--69}

     \pgap[\bibgap]

     \bibitem{La-vF2}
     \bibart{M. L. Lapidus and M. van Frankenhuijsen}
        {Fractality, Self-Similarity and Complex Dimensions}
        {Proc. Symposia Pure Math.,}
        {Part 1,}{72}{2004}{349--372}

     \pgap[\bibgap]

     \bibitem{FGCD}
     \bibbook{M. L. Lapidus and M. van Frankenhuijsen}
        {Fractal Geometry, Complex Dimensions and Zeta Functions: Geometry and spectra of fractal strings}
        {Springer Mathematical Monographs, Springer-Verlag}
        {New York}
        {2006. (2nd revised and enlarged edition to appear in 2010)}

     \pgap[\bibgap]

     \bibitem{Mat}
     \bibbook{P. Mattila}
        {Geometry of Sets and Measures in Euclidean Spaces (Fractals and Rectifiability)}
        {Cambridge Univ. Press}
        {Cambridge}
        {1995}

     \pgap[\bibgap]

     \bibitem{SST}
     \bibart{E. P. J. Pearse}
        {Canonical self-affine tilings by iterated function systems}
        {Indiana Univ. Math J.}
        {No. 6,}{56}{2007}{3151--3169. \arxiv{math/0606111}}

     \pgap[\bibgap]

     \bibitem{Pe2}
     \bibbook{E. P. J. Pearse}
        {Complex Dimensions of Self-Similar Systems}
        {Ph.D. Dissertation}
        {Univ. of California, Riverside}
        {June 2006}

     \pgap[\bibgap]

     \bibitem{GeometryOfSST}
     E. P. J. Pearse and S. Winter,
        {Geometry of self-similar tilings,}
        {to appear: Rocky Mountain J. Math. \arxiv{0811.2187}.}

     \pgap[\bibgap]

     \bibitem{Rud3}
     \bibbook{W. Rudin}
        {Functional Analysis \emph{(2nd ed.)}}
        {McGraw-Hill}
        {New York}
        {1991}

     \pgap[\bibgap]

     \bibitem{Schn1}
     \bibart{R. Schneider}
        {Curvature measures of convex bodies}
        {Ann. Mat. Pura Appl.}
        {IV}{116}{1978}{101--134}
     \pgap[\bibgap]

     \bibitem{Schn2}
     \bibbook{R. Schneider}
        {Convex Bodies: The Brunn--Minkowski Theory}
        {Cambridge Univ. Press}
        {Cambridge}
        {1993}

     \pgap[\bibgap]

     \bibitem{Stacho}
     \bibart{L. L. Stacho}
        {On curvature measures}
        {Acta Sci. Math.}
        {}{41}{1979}{191--207}

     \pgap[\bibgap]

     \bibitem{Tr}
     \bibbook{C. Tricot}
        {Curves and Fractal Dimensions}
        {Springer-Verlag}
        {New York}
        {1995}

     \pgap[\bibgap]

     \bibitem{We}
     \bibart{H. Weyl}
        {On the volume of tubes}
        {Amer. J. Math.}
        {}{61}{1939}{461--472}

     \pgap[\bibgap]

     \bibitem{We2}
     \bibart{H. Weyl}
        {\"{U}ber die Abh\={a}ngigkeit der Eigenschwingungen einer Membran von deren Begrenzung}
        {J. Reine Angew. Math.}
        {}{141}{1912}{1--11}

     \pgap[\bibgap]

     \bibitem{Za1}
     \bibart{M. Z\"{a}hle}
        {Integral and current representation of Federer's curvature measures}
        {Arch. Math.}
        {}{46}{1986}{557--567}

     \pgap[\bibgap]

     \bibitem{Za2}
     \bibart{M. Z\"{a}hle}
        {Curvatures and currents for unions of sets with positive reach}
        {Geom. Dedicata}
        {}{23}{1987}{155--171}

  \end{thebibliography}
\end{document}